%% file: MDSBlowupWPS.tex
\documentclass[11pt,letterpaper]{amsart}

\usepackage{enumitem}
\usepackage{amsfonts}
\usepackage{amssymb}
\usepackage{amsmath}
\usepackage{amsthm}
\usepackage{mathtools}
\usepackage{amstext}
\usepackage{tikz-cd}
\usepackage{enumitem}
\usepackage[left=1.3in,right=1.3in,top=1.26in, bottom=1.3in]{geometry}

\usepackage{graphicx}

\definecolor{cit}{HTML}{117733}
\definecolor{lin}{HTML}{0003d6}
\definecolor{gree}{HTML}{6dc066}
\usepackage{hyperref}
\hypersetup{bookmarksnumbered,colorlinks=true,linkcolor=lin,citecolor=cit,urlcolor=black}

\usepackage[section]{placeins}
\usepackage{caption}

\DeclareMathOperator{\bl}{Bl}

\DeclareMathOperator{\Nef}{Nef}

\DeclareMathOperator{\Pic}{Pic}
\DeclareMathOperator{\NE}{NE}
\DeclareMathOperator{\NEb}{\overline{\NE}}
\DeclareMathOperator{\Cl}{Cl}

\DeclareMathOperator{\spec}{Spec}
\DeclareMathOperator{\Hom}{Hom}

\DeclareMathOperator{\nn}{N}

\newcommand{\CC}{\mathbb{C}}
\newcommand{\Z}{\mathbb{Z}}

\newcommand{\R}{\mathbb{R}}
\newcommand{\PP}{\mathbb{P}}
\newcommand{\NN}{\mathbb{N}}

\newcommand{\Q}{\mathbb{Q}}

\newcommand{\e}{\equiv}
\newcommand{\m}{\pmod} 

\newcommand{\gkt}{of Gonz\'{a}lez-Karu type}

\newcommand{\vv}[1]{\lvert#1\rvert}
\newcommand{\vr}[1]{\langle#1\rangle}
\newcommand{\ra}{\to}
\newcommand{\xra}{\xhookrightarrow}

\newcommand{\mt}{\mapsto}

\newcommand{\du}{\vee}

\newcommand{\ik}{^{-1}}

\newcommand{\vem}{\vspace{2em}}

\newcommand{\dis}{\displaystyle}
\newcommand{\txt}[1]{\text{ #1}} 
\theoremstyle{plain}
\newtheorem{theorem}{Theorem}[section]
\newtheorem{corollary}[theorem]{Corollary}
\newtheorem{prop}[theorem]{Proposition}
\newtheorem{lemma}[theorem]{Lemma}

\theoremstyle{definition}
\newtheorem{defi}[theorem]{Definition}

\newtheorem{example}[theorem]{Example}
\newtheorem{remark}[theorem]{Remark}
\newtheorem{notation}[theorem]{Notation}
\newcommand{\pf}{{\em Proof}}
\newcommand{\pfof}[1]{{\em Proof of #1}}

\begin{document}
\title{Mori Dream Spaces and blow-ups of weighted projective spaces}
\author{Zhuang HE}
\date{}

\begin{abstract}
	For every $n\geq 3$, we find a sufficient condition for the blow-up of a weighted projective space $\PP(a,b,c,d_1,\cdots,d_{n-2})$ at the identity point not to be a Mori Dream Space.
We exhibit several infinite sequences of weights satisfying this condition in all dimensions $n\geq 3$.\end{abstract}

\maketitle
\section{Introduction}
\label{mainsection}
We study the question whether the blow-up of a projective, $\Q$-factorial toric variety over $\CC$ of Picard number one, at the identity point $p$ of the open torus, is a Mori Dream Space (MDS).

Mori Dream Spaces were introduced by Hu and Keel in \cite{hu2000}. By \cite{BCHM}, log Fano varieties over $\CC$ are Mori Dream Spaces. Projective, $\Q$-factorial toric varieties, being log Fano, are MDS. The property of being a MDS is nevertheless not a birational invariant. In fact, the blow-up of $\PP^n$ at $r$ very general points stops being a MDS if $r>8$ for $\PP^2$ and $\PP^4$, $r>7$ for $\PP^3$, and $r>n+3$ for $n\geq 5$ \cite{Mukai05}. 
One of the motivations to study blow-ups of toric varieties at the identity point comes from the proof by Castravet and Tevelev \cite{castravettevelev2015} that the moduli spaces of stable rational curves 
$\overline{M}_{0,n}$ are not MDS when $n> 133$, which was later improved to $n>12$ by Gonz\'{a}lez and Karu \cite{GK} and to $n>9$ by Hausen, Keicher and Laface \cite{hkl}. The proof of \cite{castravettevelev2015} used the examples of not MDS blow-ups of weighted projective planes (see \ref{gnwS1} and \ref{gnwS2}) by Goto, Nishida and Watanabe \cite{gnw1994}.

The discussion above prompts the question of searching for not MDS blow-ups of toric varieties of small Picard numbers, which was formulated in \cite{castravet16}.  
Historically, much research work was done for surfaces. For a weighted projective plane $S=\PP(a,b,c)$, let $p$ be the identity point of the open torus. If the anticanonical divisor $-K$ of the blow-up $\bl_p S$ of $S$ at $p$  is big, then $\bl_p S$ is a MDS \cite{cutkosky1991}. If one of $a,b,c$ is at most $4$ or equals $6$ then $\bl_p S$ is a MDS \cite{cutkosky1991}\cite{HS}. The first examples where $\bl_p S$ is not a MDS were given in \cite{gnw1994}. A generalization was achieved by Gonz\'{a}lez and Karu \cite{GK} for toric varieties of Picard number one whose corresponding polytope $\Delta$ has specific numbers of lattice points in its columns. The question can be formulated as an interpolation problem on the lattice points in $\Delta$ and leads to $3$ families of new nonexamples \cite{zh2017n}. We note that for any weighted projective space $X$, $\bl_p X$ is a MDS if and only if the Cox ring of $\bl_p X$ is a finitely generated $\CC$-algebra, which is also equivalent to the finite generation of the symbolic Rees algebra associated to $X$ \cite{cutkosky1991}\cite{gnw1994}, which is of independent interest.

In higher dimensions not much was known until the recent work \cite{GK17}. In \cite{GK17}  Gonz\'{a}lez and Karu constructed higher dimensional toric varieties $X$ of Picard number one with $\bl_p X$ not a MDS, by exhibiting a nef but not semiample divisor on $\bl_p X$. Their examples include some weighted projective $3$-spaces $X=\PP(a,b,c,d)$ such that $\bl_p X$ is not a MDS. 

\vem

In this paper, we give a sufficient condition (Theorem \ref{main}) so that the blow-up of the weighted projective $n$-space $X=\PP(a,b,c,d_1,d_2,\cdots,d_{n-2})$ at the identity $p$ is not a MDS. We show new examples of such $X$ in all dimensions $n\geq 3$.

We sum up our results below.
We work over the complex numbers $\CC$.
Let $N=\Z^2$ and $M$ be the dual lattice of $N$.
Let $S$ be a normal projective, $\Q$-factorial toric surface of Picard number $1$, with fan $\Sigma_S$ in $N\otimes_\Z \R=\R^2$. Then a polarization $H=H_\Delta$ on $S$ is determined by a rational triangle $\Delta$ in $M\otimes_\Z \R$ whose normal fan is $\Sigma_S$. Let the sides of $\Delta$ have rational slopes $s_1<s_2<s_3$. We choose $\Delta$ so that after translating one vertex of $\Delta$ to $(0,0)$, the opposite side passes through $(0,1)$. Then the {width} of this $\Delta$ equals $w:=1/(s_2-s_1)+1/(s_3-s_2)$. This $w$ is called the {\em width} of the polarized toric surface $(S,H_\Delta)$  (see \cite[Thm 1.2]{GK}).

A weighted projective plane $S=\PP(a,b,c)$ is an example of normal $\Q$-factorial toric surfaces of Picard number $1$. A triple $(e,f,-g)$ is called a {\em relation} between the weights $(a,b,c)$ if $e,f,g\in \Z_{>0}$ and $ae+bf=cg$ \cite[Thm. 1.5]{GK}. Then there exists a polarization $H_\Delta$ such that the width $w$ of $(S,H_\Delta)$ is smaller than $1$ if and only if there exists a relation $(e,f,-g)$ with $cg^2/ab=w<1$. Such $(e,f,-g)$ is unique if it exists, even when permuting the weights $a,b,c$. Therefore for a relation $(e,f,-g)$ we define the {\em width} of $(e,f,-g)$ to be $cg^2/(ab)$.

\begin{figure}
		\centering
		\includegraphics[width=8cm]{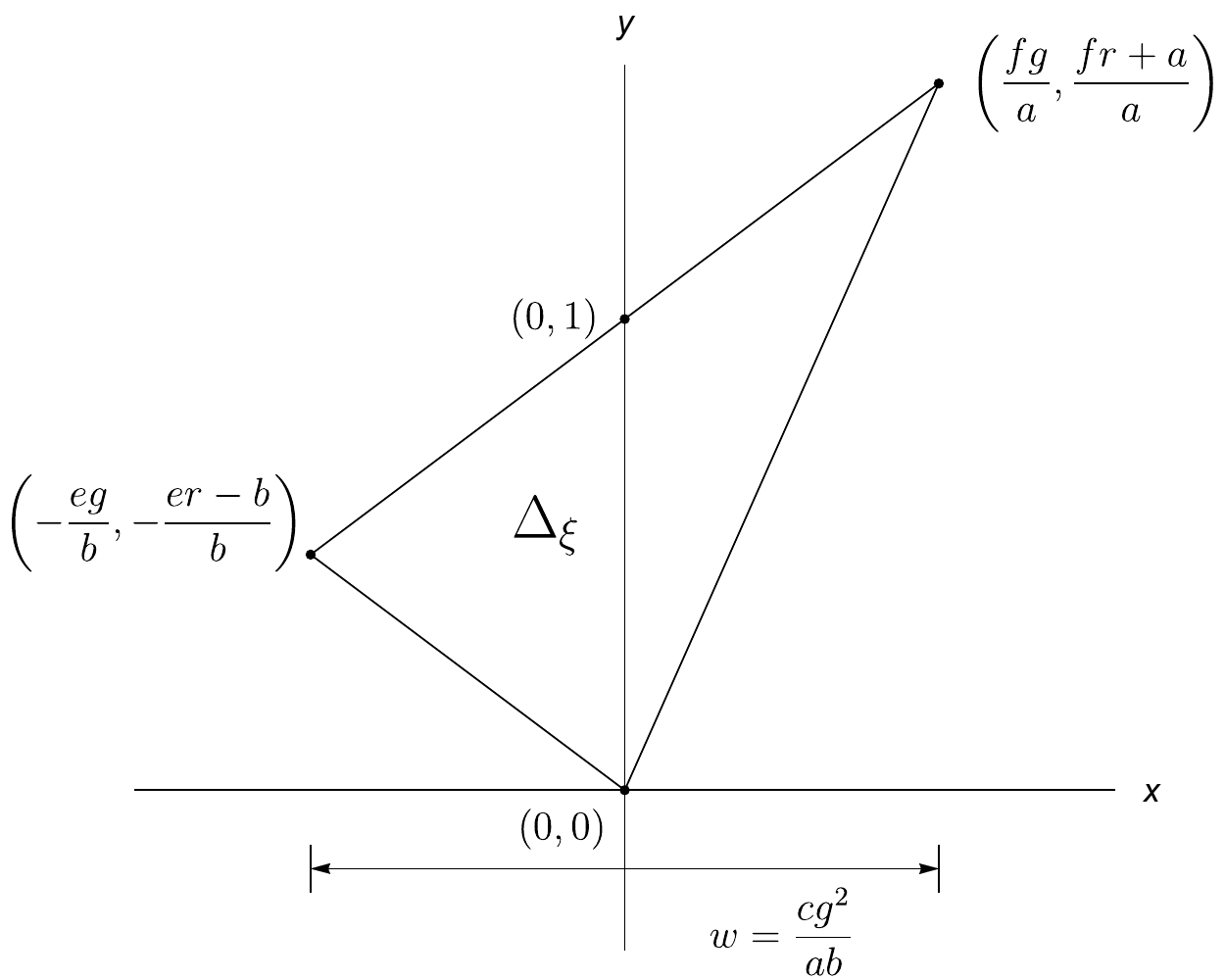}
		\caption{The triangle $\Delta_\xi$ from the relation $\xi=(e,f,-g)$}
		\label{fig:tri}
	\end{figure}

Given $\xi=(e,f,-g)$ a relation with width $w<1$, we can construct a fan $\Sigma_\xi$ of $S$ and the polytope $\Delta_\xi$ with width $w$ as follows:
By \cite[Prop. 5.1]{zh2017n}, there exists a unique integer $r$ such that $1\leq r\leq g$, $g\mid er-b$ and $g\mid fr+a$. Then the following vectors are primitive and span $\Z^2$:
	\begin{align}
		u_0=\left(\frac{er-b}{g}, -e\right), \quad u_1=\left(\frac{fr+a}{g},-f\right),\quad u_2=(-r,g).\label{Sfan}
	\end{align}
	Clearly $au_0+bu_1+cu_2=0$. Hence the fan $\Sigma_\xi$ with ray generators $u_0,u_1$ and $u_2$ is a fan of $\PP(a,b,c)$. The triangle $\Delta_\xi$ has vertices
	\begin{align}
		(0,0), \quad \left(-\frac{eg}{b},-\frac{er-b}{b}\right), \quad \left(\frac{fg}{a},\frac{fr+a}{a}\right),
		\label{Spoly}
	\end{align}
	which is normal to $\Sigma_\xi$ and has width $w=cg^2/(ab)$ (See Figure \ref{fig:tri}).

Throughout this paper, we always assume that the weights $q_0,q_1,\cdots,q_{n}$ of a weighted projective $n$-space $\PP(q_0,q_1,\cdots,q_{n})$ are well-formed, i.e., any $n$ weights are relatively prime.

For any weighted projective space $X$, let $p$ be the identity point of the open torus in $X$. For $S=\PP(a,b,c)$, let $B$ be the pseudo-effective divisor on $S$ generating $\Cl(S)\cong\Z$. Let $e$ be the exceptional divisor of the blow-up $\pi:\bl_p S\ra S$.
Our main result is:

\begin{theorem}\label{g}
	Let $X=\PP(a,b,c,d_1,d_2,\cdots,d_{n-2})$ where $a,b,c$ are pairwise coprime. Let $S=\PP(a,b,c)$. Suppose there is a negative curve $C$ on $\bl_p S$, different from $e$, with $C\sim_\Q \lambda\pi^*B-\mu e$ for some $\lambda, \mu\in\Q$. Suppose all the following hold:
	\begin{enumerate}[label={\rm(\roman*)}]
		\item every $d_i$ lies in the semigroup generated by $a,b$ and $c$ (i.e., $d_i$ is a linear combination of $a,b,c$ with non-negative integer coefficients),
		\item $\dis d_i<\frac{abc\mu}{\lambda}$ for every $i$,
		\item $\bl_p \PP(a,b,c)$ is not a MDS.
	\end{enumerate}
	Then $\bl_p X$ is not a MDS.
\end{theorem}

We show a special case of Theorem \ref{g} when there is a relation $(e,f,-g)$ between the weights $(a,b,c)$ with $w<1$. In this case, there exists a negative curve $C\sim cg\pi^*B-e$ on $\bl_p S$, and we have:

\begin{theorem}
	\label{main}
	Let $X=\PP(a,b,c,d_1,d_2,\cdots,d_{n-2})$ be a weighted projective $n$-space where $a,b,c$ are pairwise coprime. Let $p$ be the identity point of the open torus in $X$. Suppose all the following hold:
	\begin{enumerate}[label={\rm(\roman*)}]
		\item there is a relation between the weights $(a,b,c)$ such that the width satisfies $w<1$.
		\item every $d_i$ lies in the semigroup generated by $a,b$ and $c$.
		\item $d_i^2 w<abc$ for every $i$,
		\item $\bl_p \PP(a,b,c)$ is not a MDS.
	\end{enumerate}
	Then $\bl_p X$ is not a MDS.
\end{theorem}
In particular, if all $d_i=a$ and $a<b<c$ with $w<1$, then $d_i^2 w=a^2w<a^2<abc$. Thus we have the following corollary:
\begin{corollary}\label{c1}
	Assume that $a<b<c$ are pairwise coprime. Suppose $\bl_p \PP(a,b,c)$ is not a MDS, and there is a relation between the weights $(a,b,c)$ such that the width satisfies $w<1$.
Then $\bl_p \PP(a,b,c,a,\cdots,a)$ is not a MDS.
\end{corollary}

\begin{example}\label{gnwS1}
	By \cite{gnw1994}, the Cox ring of the blow-up of $\PP(a,b,c)$ at the identity point is not finitely generated as a $\CC$-algebra when $(a,b,c)=(7m-3,8m-3, (5m-2)m)$ for $m\geq 4$ and $3\nmid m$. Equivalently, the blow-up at $p$ is not a MDS. The sequence of weights has relation $(e,f,-g)=(m,m,-3)$ so that $w<1$.

	By Theorem \ref{main}, we conclude that $\bl_p \PP(7m-3,8m-3, (5m-2)m, d_1,\cdots,d_{n-2})$ is not a MDS when
	\begin{enumerate}[label={\rm(\roman*)}]
		\item $m\geq 4$ and $3\nmid m$,
		\item every $d_i$ lies in the semigroup generated by $7m-3,8m-3$ and $(5m-2)m$, and
		\item every $d_i< {(7m-3)(8m-3)}/{3}$.
	\end{enumerate}

	By Corollary \ref{c1}, $\bl_p \PP(7m-3,8m-3, (5m-2)m,7m-3,\cdots,7m-3)$ is not a MDS  for $m\geq 4$ and $3\nmid m$.

\end{example}
\begin{example}\label{gnwS2}
Another infinite sequence given by \cite{gnw1994} where the blow-ups at $p$ are not MDS is $(a,b,c)=(7m-10,8m-3, 5m^2-7m+1)$ for any $m\geq 5$ such that $3\nmid 7m-10$ and $m\not\e -7 \m{59}$ (By \cite{GK} the blow-up at $p$ is also not a MDS when $m=3$). The sequence of weights has relation $(e,f,-g)=(m,m-1,-3)$ so that $w<1$.

We conclude by Theorem \ref{main} that $\bl_p \PP(7m-10,8m-3, 5m^2-7m+1, d_1,\cdots,d_{n-2})$ is not a MDS when
	\begin{enumerate}[label={\rm(\roman*)}]
		\item $m\geq 3$, $3\nmid 7m-10$ and $m\not\e -7 \m{59}$,
		\item every $d_i$ lies in the semigroup generated by $7m-10,8m-3$ and $ 5m^2-7m+1$, and
		\item every $d_i< {(7m-10)(8m-3)}/{3}$.
	\end{enumerate}

\end{example}

\begin{example}\label{gkS3}
	The infinite sequence $(a,b,c)=(7,15+2t,26+3t)$ for $t\geq 0$ has the relation $(e,f,-g)=(1,3,-2)$. The weights $(a,b,c)$ are pairwise coprime if and only if $7\nmid t-3$. They all satisfy the criterion of \cite[Thm. 1.5]{GK}, so $\bl_p \PP(a,b,c)$ is not MDS for every $t\geq 0$, where the width
	\[w=\frac{4(26+3t)}{7(15+2t)}=\frac{104+12t}{105+14t}<1\]
	for $t\geq 0$. Theorem \ref{main} (3) then gives the upper bound
	\[d<\sqrt{\frac{abc}{w}}=\frac{ab}{g}=\frac{7(15+2t)}{2}.\]
Note that when $t\geq 0$, $a+b=\dis 2t+22<\frac{7(15+2t)}{2}$. Hence $d=a+b$ is on the list. 
As a result, $\bl_p \PP(7,15+2t,26+3t, d_1,\cdots,d_{n-2})$ is not a MDS when
	\begin{enumerate}[label={\rm(\roman*)}]
		\item $t\geq 0$ and $7\nmid t-3$,
		\item every $d_i$ lies in the semigroup generated by $7,15+2t$ and $26+3t$, and
		\item every $d_i< 7(15+2t)/2$.
	\end{enumerate}

\end{example}

The paper is organized as follows. In Section \ref{2}, we give a sufficient condition (Theorem \ref{general}) for the blow-up $\bl_p X$ of a normal projective variety $X$ with Picard number $1$ not to be a MDS, with $p$ a smooth point on $X$. Such $\bl_p X$ has a nef but not semiample divisor. 
Sections \ref{fad} and \ref{normality} consider weighted projective $n$-spaces $X$ with properties described in Theorem \ref{g}. We show that $X$ contains a closed subvariety isomorphic to $S=\PP(a,b,c)$.
Section \ref{intersection} verifies the conditions in Theorem \ref{general} for $X$ and $S$, applying a result of Fulton and Sturmfels \cite[Lem. 3.4]{FS1997}. In particular, we prove that $\bl_p X$ is not a MDS.

In Section \ref{pair}, we compare our results with the examples in \cite{GK17}. Proposition \ref{comparison} describes the overlap of our list in dimension $3$ with Gonz\'{a}lez and Karu's in \cite{GK17}. The only common examples are $X=\PP(a,b,c,cg)$ where $(e,f,-g)$ is a relation between $(a,b,c)$, and $\bl_p \PP(a,b,c)$ is not a MDS and satisfies the assumptions in \cite[Cor. 2.5]{GK17}. Note that we give more examples beyond the overlap (Examples \ref{gnwS1}, \ref{gnwS2} and \ref{gkS3}).

In Section \ref{extra}, we apply Theorem \ref{g} to the case when $X=\PP(a,b,c,d_1,d_2,\cdots,d_{n-2})$ where $S=\PP(a,b,c)$ being of the form considered in \cite[Ex. 1.4]{AGK17}. Hence $\bl_p S$ is again not a MDS. This leads to new examples where $\bl_p X$ is not MDS in Corollary \ref{mgeq1}.

\section{Blow-ups of varieties of Picard number one}\label{2}
	Let $X$ be a normal, projective, $\Q$-factorial variety of Picard number $1$ and dimension $n\geq 3$.	
	Suppose $Y_1,\cdots,Y_{n-2}$ are prime Weil divisors of $X$ $(Y_i$ not necessary normal$)$, such that the set-theoretic intersection $S:=\cap_{i=1}^{n-2} Y_i$, with the reduced subscheme structure on $S$, is a normal, projective, $\Q$-factorial surface of Picard number $1$. In addition, suppose both $\Pic(X)$ and $\Pic(S)$ are finitely generated.

	Let us blow up $S$ and $X$ at a point $p\in S$ which is smooth in $X$, $S$ and each $Y_j$. 
Let $f:\bl_p S\ra \bl_p X$ be the natural inclusion. 
Let $E$ be the exceptional divisor of the blow-up $\pi_X:\bl_p X\ra X$ and $e$ be the exceptional divisor of $\pi:\bl_p S\ra S$.

	\begin{theorem} \label{general}
		Let $X,Y_i,S$ and $f$ be defined as above.
	Suppose there exists an irreducible curve $C$ in $\bl_p S$, different from the exceptional divisor $e$ in $\bl_p S$, with $C^2<0$, such that for every $i$, $(f_* C).\bl_p Y_i<0$ in $\bl_p X$.
	Then if $\bl_p S$ is not a Mori Dream Space $($MDS$)$, then $\bl_p X$ is not a Mori Dream Space.\label{ndimension}
\end{theorem}

\pf. 
Here both $\bl_p S$ and $\bl_p X$ have Picard number $2$. Since $C^2<0$ in $\bl_p S$, $C$ spans an extremal ray of the Mori cone $\NEb(\bl_p S)$ \cite[Lem. 1.22]{kollar2008}. Since $e$ is numerically equivalent to a general line in the exceptional divisor $E$ of $\bl_p X$, $[e]$ spans an extremal ray in both $\NEb(\bl_p X)$ and $\NEb(\bl_p S)$.

Let $C_1$ be the image of $C$ in $\bl_p X$, and $e_1$ be the image of $e$ in $\bl_p X$. We show that $[C_1]$ spans the other extremal ray of $\NEb(\bl_p X)$.
Since $C$ is irreducible, $C_1$ is irreducible. Suppose towards a contradiction that $C_1$ is not extremal in $\NEb(\bl_p X)$. Then $C_1\equiv r_1 F_1+s_1 e_1$ for some effective curve $F_1$ and some rational numbers $r_1,s_1>0$.
Then there exists an irreducible component $F_2$ of $F_1$ such that $F_1\equiv r_2 F_2+ s_2 e_1$ for some rational numbers $r_2>0$ and $s_2\geq0$. Therefore we can assume at the beginning that $F_1$ is irreducible.
By assumption, $C_1\cdot \bl_p Y_i<0$ for every $i$. Since $\bl_p Y_i$ is isomorphic to the proper transform of $Y_i$ in $X$, and the class of $e_1$ is the class of a line in $E$, we have $e_1\cdot \bl_p Y_i\geq 0$. Therefore $F_1\cdot \bl_p Y_i<0$. The irreducibility assumption of $F_1$ implies that $F_1\subset \bl_p Y_i$. Run this for every $i$, and we have $F_1\subset \cap_i \bl_p Y_i= \bl_p S$.
Consider the pushforward $f_*: \nn_1(\bl_p S)\ra \nn_1(\bl_p X)$ and the pullback $f^*: \nn^1(\bl_p X)\ra \nn^1(\bl_p S)$.
Since $\nn^1(\bl_p S)$ is spanned by $[f^*H]$ and $[e]$ where $H=\pi_X^* H_0$ is the total transform of a very ample divisor $H_0$ on $X$, and $e\equiv f^*E$, we have $f^*$ is surjective.  
The dual paring between $\nn^1(\bl_p X)$ and $\nn_1(\bl_p X)$ (respectively $\nn^1(\bl_p S)$ and $\nn_1(\bl_p S)$)  is perfect. Hence $f_*$ is injective by the projection formula.
Now $f_*(C-r_1 F_1-s_1 e)\equiv C_1-r_1 F_1-s_1 e_1\equiv 0$. By injectivity, 
$C-r_1 F_1-s_1 e\equiv 0$. Then the ray  $\R_{\geq 0} [C]$ is not extremal in $\NEb(\bl_p S)$, which is a contradiction. Hence the ray $\R_{\geq0} [C_1]$ is extremal in $\NEb(\bl_p X)$.

Finally, suppose $\bl_p X$ is a MDS. Since $X$ is $\Q$-factorial, and $p$ is smooth in $X$. $\bl_p X$ is also $\Q$-factorial. Then the nef cone of $\bl_p X$ is generated by semiample divisors. In particular, there is a semiample divisor $D$ such that $D.C_1=0$. Therefore $f^*D\cdot C=f_*(f^*D\cdot C)=D\cdot f_*C=D\cdot C_1=0$ by projection formula. Hence $[f^*D]$ spans an extremal ray of $\Nef(\bl_p S)$. Now $f^*D$ is also semiample. This shows that $\bl_p S$ is a MDS.\qed

\input{fananddivisors.tex}

\input{normality.tex}

\input{intersection.tex}

\input{comparison.tex}

\input{application.tex}

\section{Acknowledgement} The author would like to express his gratitude to his advisor Ana-Maria Castravet for the continuous guidance and support, and numerous inspiring conversations. The author thanks Jos\'{e} Gonz\'{a}lez, Kalle Karu and Sam Payne for helpful discussions and suggestions. The author thanks the anonymous referees for their valuable input.

\bibliographystyle{halpha}
\bibliography{bibfile}

\end{document}

%% file: fananddivisors.tex
\section{Divisors on weighted projective spaces}\label{fad}

	In this section we construct the fan of the weighted projective $n$-space $X=\linebreak \PP(a,b,c,d_1,\cdots,d_{n-2})$ and define $n-2$ divisors $Y_j$ on $X$ for $j=3,4,\cdots,n$, under the assumption (i) of Theorem \ref{g}. Then we show that the set-theoretic intersection of those $Y_j$ equals the Zariski closure of a $2$-dimensional subtorus in $X$.

	\begin{notation}\label{allnote}
		We list some notations and terminology for later use.
		\begin{itemize}[noitemsep,nolistsep,leftmargin=*,label={$\cdot$}]
		\item For any integer $n\geq 3$, let $J:=\{3,4,\cdots,n\}$.
		\item Let $N\cong \Z^n$ ($n\geq 3$) be a lattice. Let $T_N=N\otimes_\Z \CC^*$. Then $T_N$ is a torus of dimension $n$.  Let $M=\Hom (N,\Z)$ be the dual lattice of $N$. Then $M=\Hom(T_N,\mathbb{G}_m)$, so each $u\in M$ defines a character $\chi^u$ on $T_N$.
		\item If $e_1,e_2,\cdots,e_n$ form a basis of $N$, then $e_1^*, \cdots, e_n^*$ form the dual basis of $M$. Write $\chi_j:=\chi^{e_j^*}$. Then $T_N=\spec\CC[\chi_1,\chi\ik_1,\cdots,\chi_n,\chi\ik_n]$.
	\item For any lattice $L$, define $L_\R:=L\otimes_\Z \R$.
	\item Let $N_{1}:=\Z\{e_1\}$ be the sublattice of $N$ spanned by $e_1$. Let $N_{12}:=\Z\{e_1,e_2\}$ be the sublattice spanned by $e_1$ and $e_2$. Let  $T_{1}:=N_{1}\otimes_\Z \CC^*$ and $T_{12}:=N_{12}\otimes_\Z \CC^*$ be the corresponding subtori of $T_N$. Let $M_{12}:=\Hom (N_{12},\Z)$.
	\item Let $L_j:=\Z\{e_1,e_2,\cdots,\widehat{e_j},\cdots,e_n\}$ for $j\in J$. Let $T_j:=L_j\otimes_\Z \CC^*$.
\end{itemize}

	\begin{itemize}[noitemsep,nolistsep,leftmargin=*,label={$\cdot$}]
		\item Let $\Sigma$ be a full dimensional fan in $N_\R$. If $X$ is the toric variety corresponding to the fan $\Sigma$, then $T_N$ is the open torus in $X$.
			For any full dimensional cone $\sigma\in \Sigma$, let $U_\sigma:=\spec\CC [\sigma^\du \cap M]$. Then $\{U_\sigma\mid \sigma\in \Sigma \txt{is full dimensional}\}$ is an affine open cover of $X$. 
		\item Write $\tau\prec \sigma$ if $\tau$ is a face of $\sigma$. For any cone $\tau\in\Sigma$, let $O(\tau)$ be the $T_N$-orbit associated to $\tau$ in $X$. Then for a full dimensional cone $\sigma$ and any cone $\tau$ in $\Sigma$, $O(\tau)\subseteq U_{\sigma}$ if and only if $\tau\prec\sigma$ (see \cite[3.2.6c]{coxtoric}).
		\item Let $V(\tau)$ be the Zariski closure of $O(\tau)$ in $X$. Then $V(\tau)$ is a torus-invariant closed subvariety of $X$.
		\item A fan $\Sigma$ is simplicial if any cone $\sigma\in\Sigma$ is generated by linearly independent generators. Assume that $\Sigma$ is a simplicial fan in $\R^n$ with $n+1$ rays $R_0,R_1,\cdots,R_n$, where every $n$ of them are linearly independent. For every $I\subseteq \{0,1,\cdots,n\}$, let $\sigma_I\in \Sigma$ be the cone spanned by $\{R_i\mid i\in I\}$. Every cone $\sigma\in \Sigma$ corresponds to a unique subset $I$ in the way above. Let $\Sigma(k)$ be the $k$-dimensional cones in $\Sigma$. Then $\Sigma(k)=\{\sigma_I\mid \vv{I}=k\}$. We write $V(\sigma_I)$ as $V_I$, and  $O(\sigma_I)$ as $O_I$. Then $O_I$ is a torus of dimension $n-\vv{I}$. If $I=\{i\}$, then we write the torus-invariant divisor $V(\sigma_{\{i\}})$ as $D_i$:.
	\end{itemize}
	\end{notation}

	We start with the fan of the weighted projective plane $\PP(a,b,c)$.
	The assumption and conclusion of Proposition \ref{g} are symmetric about $a,b$ and $c$. Hence up to a permutation on $(a,b,c)$, we can choose a fan $\Sigma_S$ of $S$ with ray generators $u_i=(x_i,y_i)$ such that both $y_0,y_1<0$ and $y_2>0$. Note that we have $au_0+bu_1+cu_2=0$.

	Consider $N=\Z^n$. Fix a basis $e_1,e_2,\cdots,e_n$ of $N$. By assumption (ii), there exist nonnegative integers $\{m_{ij}\}$ such that $d_{j-2}=am_{0,j} +bm_{1,j} +cm_{2,j}$ for every $j\in J$. Define the following vectors in $N$:
	\begin{align}\label{ug}
	\begin{split}
		v_0&=(x_0,y_0,-m_{0,3},\cdots,-m_{0,n}),\\
v_1&=(x_1,y_1,-m_{1,3},\cdots,-m_{1,n}),\\
v_2&=(x_2,y_2,-m_{2,3},\cdots,-m_{2,n}),\\
v_j&=e_j, \txt{for } j\in J=\{3,4,\cdots,n\}.
	\end{split}
\end{align}

Note that for every $j\in J$, at least one of the integers $m_{0j},m_{1j},m_{2j}$ is necessarily nonzero.

Those $v_i$ satisfy the relation
\[av_0+bv_1+cv_2+d_1 v_3+\cdots+d_i v_{i+2}+\cdots+d_{n-2}v_{n}=0.\]
Moreover, each $v_i$ is primitive, and together they span the lattice $N$. 
As a result, if we let $\Sigma_X$ be the fan in $N_\R$ spanned by the $n+1$ rays along $v_i$ ($i=0,1,\cdots,n$), then $\Sigma_X$ is a fan of $X=\PP(a,b,c,d_1,\cdots,d_{n-2})$.

\begin{defi}\label{YSZ}
	Let the fan $\Sigma_X$ of $X=\PP(a,b,c,d_1,\cdots,d_{n-2})$ be defined as above. For every $j\in J$, let $Y_j$ be the Zariski closure of the subtorus $T_j=L_j\otimes_\Z \CC^*$ in $X$. Define $S$ to be the set-theoretic intersection $\cap_{j=3}^n Y_j$. Let $Z$ be the Zariski closure of the subtorus $T_{12}=N_{12}\otimes_\Z \CC^*$ in $X$.  
\end{defi}

By definition, all the $Y_j$ and $Z$ are irreducible. We claim:

\begin{prop}
	\begin{enumerate}[label={\rm(\roman*)}]
		\item The set-theoretic intersection $S$ equals $Z$.
		\item With the reduced subscheme structure, $S$ is isomorphic to $\PP(a,b,c)$. In particular, $S$ is normal.
	\end{enumerate}
	\label{pureintersection}
\end{prop}

We prove (ii) of Proposition \ref{pureintersection} in the next section. Here we prove (i) by showing that $Z$ is the unique irreducible component of the intersection $S$. We will reduce the question to the affine case and apply the following lemma.
\begin{lemma}
	Let $\sigma$ in $N_\R$ be a simplicial cone spanned by $n$ linearly independent rays $R_i$, $i=1,\cdots,n$.
	Let $U_\sigma:=\spec\CC [\sigma^\du \cap M]$. For any $u\in M$ such that $u$ is primitive and $u\neq 0$, let $T_u$ be the subtorus of $T_N$ defined by $\chi^u=1$, and take the Zariski closure $\overline{T_u}$ in $U_\sigma$. Then we have:
	\begin{enumerate}[label={\rm(\roman*)}]
		\item If $\tau\prec \sigma$ such that $u\in \tau^\du \cup (- \tau^\du)$ and $u\not\in \tau^\perp$, then the set-theoretic intersection $\overline{T_u}\cap O(\tau)= \emptyset$. In particular:
	\begin{enumerate}
		\item For $\tau=R_i$, if $u\not\in \tau^\perp$, then $\overline{T_u}\cap O(\tau)= \emptyset$.
		\item If $u\in \sigma^\du \cup (- \sigma^\du)$, then $\overline{T_u}\cap O(\sigma)= \emptyset$.
	\end{enumerate}
		\item If $u\in \tau^\perp$ and $u\in \sigma^\du \cup (- \sigma^\du)$, then  $\overline{T_u}\cap O(\tau)$ has codimension at least 1 in $O(\tau)$.
	\end{enumerate}
\label{cutorbit}
\end{lemma}

\pf. When $\tau=R_i$ is a ray, $\tau^\du \cup (- \tau^\du)=M$. When $\tau=\sigma$, $\tau^\perp=\sigma^\perp=\{0\}$. Therefore the two special cases (a) and (b) of (i) follow from the general result.
Now let $\tau$ be a $d$-dimensional face of $\sigma$ such that $u\in \tau^\du \cup (- \tau^\du)$, and $u\not\in \tau^\perp$. Then $O(\tau)\cong \spec \CC[\tau^\perp \cap M]$ is a $(n-d)$-dimensional torus (see Notation \ref{allnote}).
Let $V(\tau)$ be the closure of $O(\tau)$ in $U_\sigma$. Then $V(\tau) \cong \spec\CC[\tau^\perp \cap \sigma^\du\cap M]$. Then the inclusions
\[O(\tau)\cong \spec \CC[\tau^\perp \cap M]\xra{} V(\tau) \cong \spec\CC[\tau^\perp \cap \sigma^\du\cap M]\xra{} U_\sigma\cong\spec \CC[\sigma^\du \cap M]\]
correspond to the maps of $\CC$-algebras
\[\CC[\sigma^\du\cap M]\xrightarrow{\phi_\tau} \CC[\tau^\perp \cap \sigma^\du\cap M]\ra \CC[\tau^\perp \cap M],\]
where $\phi_\tau$ sends $\chi^u$ to $\chi^u$ if $u\in \tau^\perp$, and $0$ otherwise.
To prove that $\overline{T_u}$ does not intersect $O(\tau)$, it suffices to show that there is a regular function $f$ vanishing on $T_u$ but not vanishing anywhere on $O(\tau)$.
There are two cases.

Case I.  $u\in \sigma^\du \cup (- \sigma^\du)$ and $u\not\in \tau^\perp$. Note that $\sigma^\du\subseteq \tau^\du$ since $\tau\prec \sigma$. Suppose $u\in -\sigma^\du$. Then $-u\in \sigma^\du$. By definition, $T_u=T_{-u}$, so we can assume $u\in \sigma^\du$. Now $f:=\chi^u-1=\chi^u-\chi^0\in \CC[\sigma^\du\cap M]$ is a regular function on $U_\sigma$. Since $u\not\in \tau^\perp$, $\phi_\tau (\chi^u)=0$. Since $0\in \tau^\perp$, $\phi_\tau(\chi^0)=1$. Therefore $\phi_\tau(f)=-1$ is a regular function on $V(\tau)$ which does not vanish on $O(\tau)$.

Case II. $\tau\neq \sigma$ is a proper face, $u\in \tau^\du \cup (- \tau^\du)$ and $u\not\in \sigma^\du \cup (-\sigma^\du)$ and $u\not\in \tau^\perp$. For each $i=1,\cdots,n$, let $r_i$ be the ray generator of the ray $R_i$. Without loss of generality, we can assume $\tau$ is the face spanned by $r_1,\cdots,r_d$, with $d<n$, and $u\in \tau^\du$.
Let $\vr{\cdot,\cdot}:N\times M\ra \Z$ be the dual pairing.
Then $\vr{r_i,u}\geq 0$ for $i=1,\dots,d$, with $\vr{r_i,u}>0$ for some $i\leq d$, and $\vr{r_j,u}<0$ for some $j\in \{d+1,\cdots,n\}$.
We claim there exist $p,q\in \sigma^\du \cap M-\{0\}$ and $k\in \Z_{>0}$ such that $ku=p-q$ and $q\in \tau^\perp$. 
Indeed, since $\sigma$ is simplicial,  $r_1,\cdots,r_n$ form a basis of $N\otimes_\Z \Q$. Let  $r_1^*,\cdots,r_n^*$ be the dual basis of $M\otimes_\Z \Q$. Then $u=u_1 r_1^*+\cdots+u_n r_n^*$ for rational numbers $u_i$, $i=1,\cdots,n$.  Define
\[p':=\sum_{u_i>0}u_i r_i^*, \txt{and }q':=-\sum_{u_i<0}u_i r_i^*.\]
Then $u=p'-q'$. Indeed both $p'$ and $q'$ are in $\sigma^\du$. Since $\vr{r_i,u}>0$ for some $i\leq d$, and $\vr{r_j,u}<0$ for some $j\in \{d+1,\cdots,n\}$, we have $p'\neq 0$ and $q'\neq 0$. Take any $k\in \Z_{>0}$ such that $kp'$ and $kq'$ are both in $M$. Let $p:=kp'$ and $q:=kq'$, then $ku=p-q$ and $p,q\in \sigma^\du \cap M-\{0\}$, which proves the claim.
Now let $f=\chi^q-\chi^p$. Then $f\in \CC[\sigma^\du \cap M]$. We have $f=\chi^q-\chi^p=-\chi^q(\chi^{ku}-1)$. Since $\chi^u-1$ divides $\chi^{ku}-1$, and $\chi^q$ has no poles on $T_u$, $f$ must vanish everywhere $T_u$. On the other hand, since $u\not\in \tau^\perp$ and $q\in \tau^\perp$, $p=ku+q\not\in \tau^\perp$. Therefore $\phi_\tau(\chi^p)=0$, and $\phi_\tau(f)=\phi_\tau(\chi^q)=\chi^q$. When restricted to $O(\tau)$, $\chi^q$ is a nonzero monomial in the coordinate functions on $O(\tau)$, therefore $\chi^q$ does not vanish anywhere on the torus $O(\tau)$. This proves (i).

By the symmetry between $u$ and $-u$, to prove (ii), we need only prove for the case when $u\in \tau^\perp\cap \sigma^\du$. In this case, $\phi_\tau(\chi^u)=\chi^u$, so $\chi^u-1$ is a regular function of $O(\tau)$. Now $\overline{T_u}$ is contained in the zero locus of $\chi^u-1$. By assumption, $u\neq 0$, so $\chi^u\neq 1$. Restricting to $O(\tau)$, $\chi^u\neq 1$ is a monomial of the coordinate functions on $O(\tau)$, so $\chi^u=1$ defines a  subtorus of codimension $1$ in $O(\tau)$. This proves (ii).
\qed

\vem
\pfof{Proposition \ref{pureintersection}  \em{(i)}.} By Definition \ref{YSZ}, $S$ is the set-theoretic intersection of $Y_j$, $j\in J$. Since each $Y_j$ has codimension one in $X$, the codimension of each irreducible component of $S$ in $X$ is at most $n-2$. For every $j\in J$, since $T_{12}\subseteq T_j$, $Z$ is contained in $Y_j$. Hence $Z$ is contained in $S$. Therefore it suffices to prove that $Z$ is the unique irreducible component of $S$ of dimension at least $2$. 

Here the fan $\Sigma_X$ is simplicial, spanned by ray generator $v_i$. By Notation \ref{allnote}, $\Sigma_X=\{\sigma_I\mid I\subseteq \{0,1,\cdots,n\}\}$. To prove that $Z$ is the unique irreducible component of $S$ of dimension at least $2$, we need only show that $S\cap O_I$ is contained in a curve for every  $1\leq \vv{I}\leq n-2$. Indeed, suppose $S\cap O_I$ is contained in a curve for every $1\leq \vv{I}\leq n-2$. Then $X\backslash T_N$ is a disjoint union of $T_N$-orbits $O_I$ for $1\leq \vv{I}\leq n-2$, with $\dim O_I=n-\vv{I}$. Therefore, if we assume there is some irreducible component $S'$ of $S$ disjoint from $Z$, then $S'$ is contained in $X\backslash T_N$, hence $\dim S'\leq 1$. This proves that $Z$ is the unique irreducible component of $S$ of dimension at least $2$.

It remains to show $S\cap O_I$ is contained in a curve for every  $1\leq \vv{I}\leq n-2$.
By Notation \ref{allnote}, $\{U_\sigma\mid \sigma \in \Sigma_X(n)\}$ is a torus-invariant open affine cover of $X$. 
For every $T_N$-orbit $O_I$ with $1\leq \vv{I}\leq n-2$, we choose some $\sigma'\in \Sigma_X(n)$ such that $\sigma_I\prec\sigma'$. Then $O_I\subseteq U_{\sigma'}$. By definition, $Y_j$ is the Zariski closure of $T_{e_j^*}$ in $X$. Indeed, $T_{e_j^*}\subseteq T_N\subseteq U_{\sigma'}$. Let $Y'_j$ be the restriction of $Y_j$ to this $U_{\sigma'}$. Then $Y'_j$ equals the Zariski closure of $T_{e_j^*}$ in $U_{\sigma'}$. We apply Lemma \ref{cutorbit} to $\sigma=\sigma'$, $\tau=\sigma_I$ and $u=e_j^*$. 
Recall (\ref{ug}) that $-m_{ij}\leq 0$ is the $j$-th entry of $v_i$ for $i=0,1,2$, $j\geq 3$. Define the following index sets: 
\begin{align*}
	I_+&:=I\cap \{0,1,2\},\\
	J_-&:=\{j\in J\backslash I\mid m_{ij}>0 \txt{for some } i\in I_+\},\\
	I_0&:=\{j\in I\cap J\mid m_{ij}=0 \txt{for all } i\in I_+\}.
\end{align*}
There are 4 possible cases: (a)  $I_+=\emptyset$; (b) $I_+\neq \emptyset$ and $J_-\neq \emptyset$; (c)  $I_+\neq \emptyset$ and $I_0\neq \emptyset$; and (d) $I_+\neq \emptyset$ and $J_-=I_0= \emptyset$.

In Cases (a) (b) and (c), we apply Lemma \ref{cutorbit} (i) to show that there exists $j\in J$ such that $Y'_j\cap O_I=\emptyset$ for some $j\in J$ and for every choice of $\sigma_I\prec \sigma'$. Hence $S\cap O_I=\emptyset$. For (d), we apply Lemma \ref{cutorbit} (ii) to show that $S\cap O_I$ is contained in a curve by choosing a specific $\sigma'$.

(a) $I_+=\emptyset$. Choose any $j\in I$. Then $e_j^*\not\in \sigma_I^\perp$ and $e_j^*\in \sigma_I^\du$. Apply Lemma \ref{cutorbit} (i) to any full dimensional $\sigma'$ such that $\sigma_I\prec \sigma'$, $\tau=\sigma_I$ and $u=e_j^*$. Then $Y'_j\cap O_I=\emptyset$.

(b) $I_+\neq \emptyset$ and $J_-\neq \emptyset$. Then choose any $j\in J_-$. We have $\vr{v_i,e_j^*}=-m_{ij}< 0$ for some $i\in I_+$, and $\vr{v_i,e_j^*}\leq 0$ for all $i\in I$. Hence $e_{j}^*\in -\sigma_I^\du$ and  $e_{j}^*\not\in -\sigma_I^\perp$. Therefore $Y'_{j}\cap O_I=\emptyset$.

(c) $I_+\neq \emptyset$ and $I_0\neq \emptyset$. Choose any $j\in I_0$. Then $\vr{v_j,e_j^*}=1>0$. If $i\in I$ and $i\neq j$, then either $i\in J$ or $i\in I_+$. If $i\in J$, then $v_i=e_i$ and $i\neq j$, so $\vr{v_i,e_j^*}=0$. If $i\in I_+$, then $\vr{v_i,e_j^*}=-m_{ij}=0$ since $j\in I_0$. Hence $e_j^*\in \sigma_I^\du$ and  $e_j^*\not\in \sigma_I^\perp$, so $Y'_j\cap O_I=\emptyset$.

(d) $I_+\neq \emptyset$ and $J_-=I_0= \emptyset$. Since $\vv{I}\leq n-2$, and $I_+\neq \emptyset$, it must be that $J\not\subseteq I$. Therefore $I_+\neq \{0,1,2\}$ (otherwise for every $j\in J\backslash I$, there exists an $m_{ij}>0$, so $j\in J_-$), so $\vv{I_+}=1$ or $2$. Fix some $j\in J\backslash I$. Since $J_-=\emptyset$, $m_{ij}=0$ for all $i\in I_+$. Therefore $e_j^*\in \sigma_I^\perp$. For this $j\in J\backslash I$, define $I'=\{0,1,2,\cdots,\widehat{j},\cdots,n\}$ and let $\sigma':=\sigma_{I'}$. Define $Y_j'$ to be the restriction of $Y_j$ to $U_{\sigma'}$ as discussed above. Then $U_{\sigma'}$ contains $O_I$, with $e_j^*\in -(\sigma')^\du$. In Lemma \ref{cutorbit} (ii), let $\sigma=\sigma'$, $\tau=\sigma_I$ and $u=e_j^*$. Then $Y'_j\cap O_I$ is of codimension at least one in $O_I$ and is contained in the zero locus of $\chi_j-1$, regarded as a regular function on $O_I$.
Now the number of such $j$ equals $\vv{J\backslash I}=n-2-\vv{I\cap J}=n-2-(\vv{I}-\vv{I_+})$. Since $n-\vv{I}=\dim O_I$, we have $ \vv{J\backslash I}=\dim O_I-(2-\vv{I_+})$. 
Recall that $M=\Z\{e_1^*,\cdots,e_n^*\}$ and $O_I=\spec \CC[\sigma_I^\perp\cap M]$. Since $\vv{I_+}=1$ or $2$, the semigroup $\sigma_I^\perp\cap M$ is generated by $\{e_i^*\mid i\in J\backslash I\}$ if $\vv{I_+}=2$, or by $\{e_i^*\mid i\in J\backslash I\}$ together with some $\xi\in \Z\{e_1^*,e_2^*\}$ if $\vv{I_+}=1$.
Therefore each $\chi_j$, $j\in J\backslash I$ restricts to different coordinate functions on $O_I$. Hence, the intersection of the zero loci of all those $\chi_j-1$ ($j\in J\backslash I$) has dimension exactly $2-\vv{I_+}$, which is either $1$ or $0$. Therefore $S\cap O_I$ is contained in a curve. This finishes Case (d) and the proof.\qed

%% file: normality.tex
\section{Normality of the closure of subtori}\label{normality}
In this section we prove (ii) of Proposition \ref{pureintersection}, namely that the surface $S$ is normal and isomorphic to the weighted projective plane $\PP(a,b,c)$.

We recall the following construction in \cite[\S 2.1]{coxtoric} of a projective toric variety $X_A$ out of a finite set of lattice points $A\subset M$.
Let $N=\Z^n$ and $M=\Hom (N,\Z)$. Then each $m\in M$ gives a character $\chi^m$ of the torus $T_N$. Any list of $k$ lattice points $A=(m_1,\cdots,m_k)\subset M$ defines a morphism $\phi_A$ from $T_N$ to $\PP^{k-1}$:
\begin{align}
	\begin{split}
		\phi_A:\quad &T_N\ra T_{k}\xrightarrow{\mu} \PP^{k-1},\\
	&t\mt (\chi^{m_1}(t),\cdots,\chi^{m_k}(t))\mt [\chi^{m_1}(t):\cdots:\chi^{m_k}(t)].
	\label{character}
	\end{split}
\end{align}
where $T_k\cong(\CC^*)^k$ and $\mu:T_{k}\ra \PP^{k-1}$ maps $T_{k}$ to the open torus $\{[x_0:\cdots:x_{k-1}]\mid\txt{all } x_i\neq 0\}$ of $\PP^{k-1}$.
\begin{defi}\cite[Definition 2.1.1]{coxtoric}\label{XA}
We denote by $X_A$ the not necessarily normal toric variety given by the Zariski closure of the image $\phi_A(T_N)$ in $\PP^{k-1}$. 
\end{defi}
\begin{remark}
	Up to isomorphism, the definition of $X_A$ only depends on the set of points appearing in $A$. So up to isomorphism we can ignore the order of the points in $A$, and can remove possible duplicates from $A$.
	\label{noorder}
\end{remark}

We note that by definition, $X_A$ is projective. However $X_A$ need not be normal. One of the ways to obtain normal toric varieties is from polytopes. Let $P$ be a full dimension polytope in $M_\R$. Call $P$ a lattice polytope if the vertices of $P$ are in $M$.
Now consider a semigroup $\mathcal{S}\subset M$, with the addition inherited from $M$. Recall that $\mathcal{S}$ is said to be {\em saturated} if for every $m\in M$, every $k\in\Z-\{0\}$, $km\in \mathcal{S}$ implies $m\in \mathcal{S}$.

\begin{defi}\cite[Definition 2.2.17]{coxtoric} A lattice polytope $M$ is very ample if for every vertex $m\in P$, the semigroup $\mathcal{S}_{P,m}$ generated by the set $P\cap M-m$ is saturated in $M$.
\end{defi}

\begin{lemma}\label{kPvample}\cite[Cor. 2.2.19]{coxtoric}
	If $P$ is a full dimensional lattice polytope, then $kP$ is very ample if $k\geq \dim P-1$. In particular, if $P$ is a lattice polygon in $\R^2$ then $P$ is very ample.
\end{lemma}

\begin{defi}\cite[Definition 2.3.14]{coxtoric}\label{XP}
	Suppose that $P\subset M_\R$ is a full dimensional lattice polytope.
	Then define the toric variety $X_P$ to be $X_A$ with $A=kP\cap M$, for any integer $k>0$ such that $kP$ is very ample.
\end{defi}
The toric variety $X_P$ is well defined since $X_{kP\cap M}$ and $X_{\ell P\cap M}$ are isomorphic when both $kP$ and $\ell P$ are very ample (see \cite[\S 2.3]{coxtoric}).

\begin{lemma}\label{vamnormal} If $P$ is a full dimensional very ample lattice polytope, then
	$X_{P\cap M}$ is a normal projective toric variety, whose fan in $N$ is the normal fan $\Sigma$ of $P$. 
\end{lemma}
\pf. This follows from \cite[Thm. 2.3.1, Thm. 1.3.5]{coxtoric}.\qed

\vem
Now we are ready to prove that $S$ is normal and isomorphic to $\PP(a,b,c)$.

\pfof{Proposition \ref{pureintersection} \em{(ii)}.} 
Let $M_{12}=\Z\{e_1^*,e_2^*\}$. We first show that $S$ is a normal projective variety. By Lemma \ref{vamnormal}, we need only show $S\cong X_{Q\cap M_{12}}$ for some full dimensional very ample lattice polytope $Q$ in $(M_{12})_\R$.
Consider $X=\PP(a,b,c,d_1,\cdots,d_{n-2})$, with the fan $\Sigma_X$ defined by generators $v_i$ in (\ref{ug}). Choose any lattice polytope $P$ in $M_\R$ whose normal fan is $\Sigma_X$. 
By replacing $P$ with some multiple $kP$, we can assume $P$ is very ample.
By Lemma \ref{vamnormal}, we have $X=X_P=X_{P\cap M}$. Let $m_0,m_1,\cdots,m_u$ be the distinct lattice points of $P\cap M$. Let $\psi:=\phi_{P\cap M}$ be the map defined in (\ref{character}). Then
\begin{align*}
	\begin{split}
\psi=\phi_{P\cap M}:\quad &T_N\ra  T_{u+1}\ra \PP^u,\\
&t\mt (\chi^{m_0}(t),\chi^{m_1}(t),\cdots,\chi^{m_u}(t)) \mt [\chi^{m_0}(t):\chi^{m_1}(t):\cdots:\chi^{m_u}(t)].
	\end{split}
\end{align*}
Then $X$ equals the Zariski closure of $\psi(T_N)$ in $\PP^u$.
Let $\rho:M\ra M_{12}$ be the projection map. If $t\in T_{12}$, then $\chi^{m_i}(t)=\chi^{\rho(m_i)}(t)$ for every $i$. Therefore, the restriction of $\psi$ on $T_{12}$ equals
\begin{align*}
	\begin{split}
		\psi_{\mid T_{12}}:\quad &T_{12}\ra T_{u+1}\ra \PP^u,\\
&t\mt (\chi^{\rho(m_0)}(t),\cdots,\chi^{\rho(m_u)}(t))\mt [\chi^{\rho(m_0)}(t):\cdots:\chi^{\rho(m_u)}(t)].
	\end{split}
\end{align*}
By Proposition \ref{pureintersection} (i), $S$ equals to the Zariski closure of $\psi(T_{12})$ in $X$. Since $X$ is closed in $\PP^u$, we have $S$ equals the Zariski closure of $\psi(T_{12})$ in $\PP^u$. 

Define $A:=\rho(P\cap M)$. Then $A$ is the set of distinct elements in the list  $A'=(\rho(m_1),\cdots,\rho(m_u))$. By Remark \ref{noorder}, we can remove the duplicates in $A'$, so that $S\cong X_A$.

Now we only need to show that $\rho(P\cap M)=\rho(P)\cap M_{12}$ and $Q:=\rho(P)$ is a full dimensional very ample lattice polytope in $M_{12}$.
We first show that $Q$ is a lattice triangle in $(M_{12})_\R$. Recall that $P$ has the following facet presentation:
\begin{equation}\label{Pfp}
	P=\{z\in M_\R\mid \vr{v_i,z}\leq a_i \txt{for } i=0,1,\cdots,n\}
\end{equation}
for some $a_i\in \Z$ (See \cite[p. 66]{introTV}, \cite[2.2.1]{coxtoric}). Since the normal fan of $P$ is $\Sigma_X$, $P$ has exactly $n+1$ facets $F_i$ whose outer normal vectors are $v_i$, $i=0,\cdots,n$ respectively.  The reason that $a_i\in \Z$ is as follows: Fix $i\in \{0,1,\cdots,n\}$. Let $m$ be a vertex of the facet $F_i$. Then $m$ is a vertex of $P$, so $m\in M$. Since $m\in F_i$, we in fact have $\vr{v_i,m}=a_i$. Thus $a_i\in \Z$ since $v_i\in N$.

Let $z=(z_1,\cdots,z_n)\in M_\R$. Then $\rho(z)=(z_1,z_2)$. By definition of $u_i$ and $v_i$ in (\ref{ug}), we have $\vr{v_i,z}=\vr{u_i,\rho(z)}-(z_3 m_{i,3}+\cdots+z_n m_{i,n})$ for $i=0,1,2$, and  $\vr{v_j,z}=z_j$ for $j\in J=\{3,4,\cdots,n\}$. Therefore $z\in P$ if and only if $\vr{u_i,\rho(z)}\leq a_i+(z_3 m_{i,3}+\cdots+z_n m_{i,n})$ for $i=0,1,2$ and $z_j\leq a_j$ for $j\in J$. Recall that every $m_{i,j}\geq 0$. As a result, $y\in Q$ if and only if $\vr{u_i,y}\leq a_i+(a_3 m_{i,3}+\cdots+a_n m_{i,n})$ for $i=0,1,2$. Define $q_i:=a_i+(a_3 m_{i,3}+\cdots+a_n m_{i,n})$ for $i=0,1,2$. Then 
\begin{equation}\label{Qfp}
	Q=\{y\in (M_{12})_\R\mid \vr{u_i,y}\leq q_i, \txt{for } i=0,1,2\}.
\end{equation}
Indeed (5) is a facet presentation of $Q$. Thus $Q$ is a triangle in $(M_{12})_\R$.

It remains to show that $Q$ is a lattice triangle. 
A point $z\in P$ (or $Q$) is a vertex of $P$ (or  $Q$) if and only if  $z$ lives in all but one facets. By the facet presentation (\ref{Pfp}) of $P$, $m$ is a vertex of $P$ if and only if $\vr{v_i,m}=a_i$ for all $v_i$ but one. 
Suppose that $\xi_0,\xi_1,\xi_2$ are the vertices of $P$ where $\xi_i$ lives in the $n$ facets except $F_i$. We claim that $\rho(\xi_0),\rho(\xi_1)$ and $\rho(\xi_2)$ are the three vertices of $Q$.
Indeed, we need only to prove this for $\xi_0$. Let $\xi_0=(z_1,\cdots,z_n)$. Then $a_j=\vr{v_j,\xi_0}=z_j$ for $j\in J$, and $a_{k}=\vr{v_{k},\xi_0}=\vr{u_{k},\rho(\xi_0)}-(z_3 m_{k,3}+\cdots+z_n m_{k,n})$ for $k=1,2$. By definition, this shows that $\vr{u_{k},\rho(\xi_0)}=q_k$ for $k=1,2$.
Let $F'_i$ be the facet of $Q$ normal to $u_i$, for $i=0,1,2$ (see (\ref{Qfp})).
Then $\rho(\xi_0)=F'_1\cap F'_2$ is a vertex of $Q$. Since $P$ is a lattice polytope, $\xi_0\in M$, so $\rho(\xi_0)\in M_{12}$. Repeat this argument for $\xi_1$ and $\xi_2$. Then $\rho(\xi_0),\rho(\xi_1)$ and $\rho(\xi_2)$ are distinct vertices of $Q$. Therefore $Q$ is a lattice triangle. 
By Lemma \ref{kPvample}, any lattice triangle in $M_{12}$ is very ample, so $Q$ is very ample. Hence we verified that $Q$ is a full dimensional very ample lattice polytope.

It remains to show $\rho(P\cap M)=\rho(P)\cap M_{12}$. By definition, $\rho(P\cap M)\subseteq \rho(P)\cap M_{12}$. Conversely, suppose $y=(z_1,z_2)\in \rho(P)\cap M_{12}$. Then $y=\rho(z)$ where $z:=(z_1,z_2,a_3,\cdots,a_n)$. By (\ref{Qfp}), we have $ \vr{u_i,y}\leq q_i$ for $i=0,1,2$. Hence $\vr{u_i,\rho(z)}\leq q_i= a_i+(a_3 m_{i,3}+\cdots+a_n m_{i,n})$ for $i=0,1,2$. The argument preceding (\ref{Qfp}) shows that $z\in P$. Since $z_1,z_2$, all $a_i$ and all $m_{i,j}$ are integers, we have $z\in M$. Thus $\rho(P)\cap M_{12}\subseteq \rho(P\cap M)$. We conclude that $\rho(P\cap M)=\rho(P)\cap M_{12}$. Therefore, $S=X_{\rho(P)\cap M_{12}}$ is normal. Furthermore, by Proposition \ref{vamnormal}, the fan of $S$ in $N_{12}$ is the normal fan of $Q$ with respect to $N_{12}$, hence is spanned by $u_0,u_1$ and $u_2$. By (\ref{ug}), the fan spanned by $u_0,u_1$ and $u_2$ is a fan of $\PP(a,b,c)$. As a conclusion, $S\cong \PP(a,b,c)$.\qed

%% file: intersection.tex
\section{Intersection products on weighted projective spaces}\label{intersection}
We prove Theorem \ref{g} and Theorem \ref{main} in this section. 
In Section \ref{fad} we constructed a fan $\Sigma_X$ for $X=\PP(a,b,c,d_1,\cdots,d_{n-2})$, under the assumption (i) of Theorem \ref{g}. Recall that $S$ is defined as the intersection of $Y_j$ for $j\in J$, where $J=\{3,4,\cdots,n\}$. By Lemma \ref{pureintersection} (ii), $S$ is isomorphic to $\PP(a,b,c)$.

We start with a review of the intersection products of various torus-invariant divisors on $X$ and $S$.
Let  $A_d(X)$ be the Chow group of $d$-dimensional cycles in $X$. Since $X$ is a complete simplicial toric variety, by \cite[Lem. 12.5.1]{coxtoric}, $A_d(X)$ is generated by the classes of torus-invariant subvarieties $[V_I]$ where $\vv{I}=n-d$. In particular, $A_{n-1}(X)$ is generated by the classes of torus-invariant Weil divisors $\{[D_i]\mid i=0,1,2,\cdots,n\}$.
The divisor class group  $\Cl(X)$ of $X$ is isomorphic to $\Z$ by \cite[Ex. 4.1.5]{coxtoric}. Let $A$ be a pseudo-effective Weil divisor on $X$ which generates $\Cl(X)$. Then in $A_{n-1}(X)=\Cl(X)$ we have
\begin{align}
	[D_0]=a[A], \quad [D_1]= b[A], \quad [D_2]= c[A], \quad [D_j]= d_{j-2}[A], \txt{for } j\geq 3.	
	\label{XB}
\end{align}
Now $\Sigma_X$ is simplicial (Notation \ref{allnote}). By \cite[Lem. 12.5.2]{coxtoric}, we have the following intersection products:
\begin{align}
	\begin{split}	
		&[A]^n=\frac{1}{abcd_1\cdots d_{n-2}},\\
		&[D_3]\cdot[D_4]\cdot\ldots\cdot[D_n]=[V_J],\\
		&[V_J]\cdot[D_i]=[V_{J\cup \{i\}}], \txt{for } i=0,1,2,\\
		&[D_1]\cdot[D_2]\cdot[V_J]=\frac{1}{a},\quad [D_0]\cdot[D_2]\cdot[V_J]=\frac{1}{b},\quad [D_0]\cdot[D_1]\cdot[V_J]=\frac{1}{c}\cdot
	\label{XI}
	\end{split}
\end{align}
By Notation \ref{allnote}, $N_{12}=\Z\{e_1,e_2\}$. Let $\Sigma_S$ in $(N_{12})_\R$ be the fan of $S$ generated by ray generators $u_0,u_1$ and $u_2$ (See (\ref{ug})). Define $B_i:=V(\sigma_{\{i\}})$ to be the torus-invariant divisors of $S$ corresponding to $u_i$. By \cite[Ex. 4.1.5]{coxtoric}, $\Cl(S)\cong\Z$. Let $B$ be a pseudo-effective Weil divisor on $S$ which generates $\Cl(S)$. Then
\begin{align}
	[B_0]= a[B], \quad [B_1]= b[B], \quad [B_2]= c[B], \quad [B]^2=\frac{1}{abc}.
	\label{SI}
\end{align}
Next we recall a result by Fulton and Sturmfels \cite{FS1997}. Let $W$ be a toric variety of a fan $\Sigma\subset N=\Z^n$. As in \cite{FS1997}, define $N_\sigma$ as $\Z(N\cap \sigma)$, the sublattice spanned by $\sigma$ in $N$.
Let $L$ be a saturated $d$-dimensional sublattice of $N$. Let $Y$ be the Zariski closure of the subtorus $T_L=L\otimes_\Z \CC^*$ in $W$.
For every lattice point $w\in N$, define
\[\Sigma(w):=\{\sigma\in \Sigma: L_\R+w \txt{meets $\sigma$ in exactly one point}\}.\]
Here $L_\R+w:=\{x+w\mid x\in L_\R\}$.
\begin{defi}\cite[\S 3]{FS1997}
	$w$ is called {\em generic} (with respect to $L$) if $\dim \sigma=n-d$ for all $\sigma\in \Sigma(w)$.
\end{defi}
\begin{lemma}
	\cite[Lem. 3.4]{FS1997} Let $W$, $L$ and $Y$ be defined as above.
	If $w\in N$ is a generic point with respect to $L$, then 
	\[[Y]=\sum_{\sigma\in \Sigma(w)} m_\sigma [V(\sigma)]\in A_d (W),\]
	where $m_\sigma:=[N:L+N_\sigma]$ is the index of the lattice sum $L+N_\sigma$ in $N$.
	\label{FS}
\end{lemma}

For simplicity, when there are no ambiguity of the choice of $L$, and when the toric variety $W$ has a simplicial fan $\Sigma$ spanned by rays $r_0,r_1,\cdots,r_n$, we write $m_{\sigma_I}=[N:L+N_\sigma]$ as $m_I$, for $I\subset\{0,1,\cdots,n\}$. When $I=\{i\}$, we write $m_{\sigma_I}$ as  $m_i$.
\begin{lemma}
	\label{Y}
	Let $X$, $Y_j$ and $S$ be defined as in Definition \ref{YSZ}. Then $[Y_j]=[D_j]$ for all $j\in J$, and $[S]=[V_J]$.
\end{lemma}
\pf. Fix $j\in J$. By Notation \ref{allnote}, $L_j:=\Z\{e_1,e_2,\cdots,\widehat{e_j},\cdots,e_n\}$. 
By Definition \ref{YSZ}, $Y_j$ is the Zariski closure of $T_j=L_j\otimes_\Z \CC^*$ in $X$. We apply Lemma \ref{FS} to $W=X$, $Y=Y_j$ and $L=L_j$. First, $e_j$ is generic with respect to $L_j$. Indeed if $j\not\in I$, then $(L_j)_\R+e_j$ does not meet $\sigma_I$. If $j\in I$, then $\sigma_I$ intersects $(L_j)_\R+e_j$ at a single point if and only if $I=\{j\}$.
Hence $\Sigma(e_j)=\{\sigma_{\{j\}}\}$. Since $\sigma_{\{j\}}$ is a $1$-dimensional cone, $e_j$ is generic. By Lemma \ref{FS}, $[Y_j]=m_j [D_j]$, and $m_j$ equals the index of $L_j+N_{\sigma_{\{j\}}}$ in $N$, which equals to $1$, so $[Y_j]=[D_j]$.

Similarly, $N_{12}:=\Z\{e_1,e_2\}$, and $S$ is the Zariski closure of $T_{12}:=N_{12}\otimes_\Z \CC^*$. The same argument above shows that $\Sigma(\omega)=\{\sigma_{J}\}$, where $\omega=(0,0,1,\cdots,1)\in N$ is generic with respect to $N_{12}$. Apply Lemma \ref{FS} to $W=X,Y=S$ and $L=N_{12}$. Then we have $[S]=m_J[V_J]$. Here $m_J=1$ since $N_{12}+N_{\sigma_J}=N$.
\qed

\begin{defi}\label{CC0}
	Let $N_1=\Z\{e_1\}$ and $T_1:=N_1\otimes_\Z \CC^*$.
	Let $C_1$ be the Zariski closure of the subtorus $T_1$ in $S$. 
\end{defi}

\begin{lemma}
	Let $C_1$ be defined as above. Then 
	\begin{enumerate}[label={\rm (\roman*)}]
		\item The irreducible curve $C_1$ equals the closure of the subtorus $T_1$ in $X$.
		\item The class $[C_1]=-y_0[V_{J\cup \{0\}}]-y_1[V_{J\cup\{1\}}]\in A_{1}(X)$.
		\item The class $[C_1]=y_2[V_J]\cdot[D_2]\in A_{1}(X)$.
		\item The class $[C_1]=y_2[B_2]=cy_2[B]\in A_{1}(S)$.
	\end{enumerate}
\label{curveclass}
\end{lemma}
\pf. Let $\overline{T_1}$ be the closure of $T_1$ in $X$.
By definition, $T_1$ is contained in $S$. Since $S$ is closed in $X$, $\overline{T_1}$ is contained in $S$. Therefore $C_1=\overline{T_1}$. Hence, both $C_1$ and $C$ are irreducible. This proves (i).
For (ii), we work in $N=\Z^n$. Define $w=(w_1,w_2,1,\cdots,1)\in N$ such that $(w_1,w_2)$ lies in the interior of the cone spanned by $u_0$ and $u_1$. We claim that $w$ is generic with respect to $N_1$.
Indeed, by the definition of $u_i$ (see (\ref{ug})), the second coordinates of $u_0$ and $u_1$ are negative and the second coordinate of $u_2$ is positive. Hence $w_2<0$. Suppose the line $\ell:=(N_1)_\R+w$ intersects $\sigma_I$. Then $J\subset I$. Since $w_2<0$, $\ell$ misses $\sigma_J$ and $\sigma_{J\cup\{2\}}$, and meets $\sigma_{J\cup\{0\}}$ and $\sigma_{J\cup\{1\}}$  at a unique point. In the remaining case, $I=J\cup \{i_1,i_2\}$ with distinct $i_1, i_2\in \{0,1,2\}$, so $\ell$ intersects $\sigma_I$ at infinitely many points. As a conclusion, $\Sigma(w)=\{\sigma_{J\cup\{0\}}, \sigma_{J\cup\{1\}}\}$, so $w$ is generic.

Apply Lemma \ref{FS} to $W=X,Y=C_1$ and $L=N_{1}$. We have 
\[[C_1]=m_{J\cup\{0\}} [V_{J\cup \{0\}}]+m_{J\cup\{1\}}  [V_{J\cup\{1\}}].\]
By definition, $m_{J\cup\{0\}}=[N:N_1+N_{\sigma_{J\cup\{0\}}}]$. Since $N_1+N_{\sigma_{J\cup\{0\}}}$ is spanned by $e_1,e_3,\cdots,e_n$ together with $v_0$, the index equals to the absolute value of the second coordinate of $v_0$, That is, $m_{J\cup\{0\}}=\vv{y_0}$. Recall our assumption in Section \ref{fad} that $y_0, y_1<0$ and $y_2>0$. Hence $m_{J\cup\{0\}}=-y_0$. Similarly we have $m_{J\cup\{1\}}=-y_1$. This proves (ii). Now use formulas (\ref{XB}) and (\ref{XI}):
\begin{align*}
[C_1]&=-y_0[V_{J\cup \{0\}}]-y_1 [V_{J\cup\{1\}}]=-y_0 [V_J]\cdot[D_0]-y_1 [V_J]\cdot[D_1]\\
&=[V_J]\cdot[-y_0 a[A]-y_1 b[A]]=cy_2 [V_J]\cdot[A]=y_2 [V_J]\cdot[D_2].
\end{align*}
This proves (iii).

Finally consider $C_1$ as a curve on $S$.
The fan $\Sigma_S$ lives in $(N_{12})_\R$ (See Notation \ref{allnote}). We have $\Sigma(e_2)=\{B_{2}\}$. Therefore $e_2=(0,1)$ is generic with respect to $N_1$. Apply Lemma \ref{FS} to $W=S,Y=C_1$ and $L=N_{1}$. Then $[C_1]=m_2[B_2]\in A_1(S)$ where $m_2=[\Z^2:(N_1)_\R+\Z u_2]=\vv{y_2}=y_2$. This proves (iv). 
\qed

\begin{lemma} Consider the class $[B]\in A_1(X)$. Then we have $\dis [B].[Y_j]=\frac{d_{j-2}}{abc}$, for $j\in J$.
	\label{bh}
\end{lemma}
\pf. By Lemma \ref{curveclass}, $ [C_1]=cy_2 [V_J]\cdot[A]\in A_1(X)$, and $[C_1]=cy_2[B]\in A_1(S)$. Therefore $cy_2 [B]= cy_2 [V_J]\cdot [A]$ in $A_1(X)$, so  $\dis [B]=[V_J]\cdot [A]=\frac{1}{a} [V_J]\cdot[D_0]$ in $A_1(X)$. Then
\[[B].[Y_j]=\frac{1}{a}[V_J]\cdot [D_0]\cdot \frac{d_{j-2}}{b}[D_1]=\frac{ d_{j-2}}{abc}.\]\qed

Now we prove Theorem \ref{g}.

\pfof{Theorem \ref{g}.} By definition $X=\PP(a,b,c,d_1,\cdots,d_{n-2})$ is a weighted projective $n$-space. By Proposition \ref{pureintersection}, $S=\PP(a,b,c)$ is a weighted projective plane. Hence both $X$ and $S$ are normal projective $\Q$-factorial varieties, with finitely generated Picard groups. 
By Proposition \ref{pureintersection}, $S=\cap_{j=3}^n Y_j$. By assumption, $C$ is a negative curve on $\bl_p S$ and $C\neq e$. To apply Theorem \ref{general} to $X, Y_j$, $S$ and $C$, we need only verify that  $(f_* C)\cdot \bl_p Y_j<0$ for $j=3,4,\cdots,n$. 
Here $ (f_* C)\cdot \bl_p Y_j=f_*C\cdot (\pi^*_X Y_j-E)$, and $C\sim_\Q \lambda \pi^*B-\mu e$. Hence by Lemma \ref{bh} and projection formula:
\begin{align*}
	f_*C\cdot (\pi^*_X Y_j-E)&=(\pi_X)_* f_*[C]\cdot [Y_j]-f_*[C]\cdot [E]\\
 	&=\lambda [B].[Y_j]-\mu=\frac{\lambda d_{j-2}}{abc}-\mu<1.
\end{align*}
By Theorem \ref{general}, $\bl_p X$ is not a MDS.
This proves the theorem.\qed

Finally we prove Theorem \ref{main}. 

\pfof{Theorem \ref{main}.} Suppose there is a relation $(e,f,-g)$ between the weights $(a,b,c)$ such that the width $w=cg^2/(ab)<1$.

We need only show that there exists a non-exceptional negative curve $C$ on $\bl_p S$ satisfying the assumption in Theorem \ref{g} with $\lambda=cg$ and $\mu=1$, and $d_i<abc\mu/\lambda=ab/g$ for all $i=0,1,\cdots n-2$.
We first choose a specific fan $\Sigma_S$ and use $\Sigma_S$ to define $\Sigma_X$.
Indeed, by \cite[Prop. 5.1]{zh2017n}, there exists a unique integer $r$ with $1\leq r\leq g$, $g\mid er-b$ and $g\mid fr+a$. Let $u_i=(x_i, y_i)$ be given by (\ref{Sfan}):
	\begin{align}
		u_0=\left(\frac{er-b}{g}, -e\right), \quad u_1=\left(\frac{fr+a}{g},-f\right),\quad u_2=(-r,g).
	\end{align}
	Then $u_i$ span a fan of $S$. 
	Let this fan be $\Sigma_S$. We check that $y_0=-e<0$, $y_1=-f<0$ and $y_2=g>0$, so all the assumptions in Section \ref{fad} are satisfied. Then we can use $u_i$ to define $v_i$ and the fan $\Sigma_X$ as in (\ref{ug}). Consider the curve $C_1$ in Definition \ref{CC0}. Let $C$ be the proper transform of $C_1$ in $\bl_p S$.
	Then $C\sim \pi^* C_1-e$ on $\bl_p S$. By Lemma \ref{curveclass} (iv), $C\sim cg \pi^*B-e$. Hence $\lambda=cg$ and $\mu=1$. By (\ref{SI}), $[B]^2=1/abc$. Hence $[C_1]^2=g^2 c^2/abc=cg^2/ab=w$, and $[C]^2=[C_1]^2-1=w-1<0$. Since $\pi(C)=C_1$ is not a point, $C$ is not $e$. As a result, $C$ is a non-exceptional negative curve on $\bl_p S$. Finally by assumption (ii) of Theorem \ref{main}, for every $i$,  $d_i^2 w<abc$. Therefore $d_i^2 cg^2/(ab)<abc$. That is, $d_i<ab/g$. By Theorem \ref{g}, we conclude that $\bl_p X$ is not a MDS.
\qed

%% file: comparison.tex
\section{Comparison with Gonz\'{a}lez and Karu's examples}\label{pair}
We compare the $3$-dimensional case of Theorem \ref{main} with  \cite[Thm. 2.3, Cor. 2.5]{GK17}. 

\begin{defi}\label{gkt}
	Consider a $n$-dimensional convex polytope $\Delta$ in $\R^n$ such that all its vertices have rational coordinates. 
	\begin{enumerate}[label={\rm(\roman*)}]
		\item For $n=3$, we say such a polytope is {\em \gkt} if the vertices of $\Delta$ are $(0,0,1)$, $(0,1,0)$, $P_L$ and $P_R$, with $P_L$ and $P_R$ and $0$ collinear, and $x(P_L)<0<x(P_R)\leq x(P_L)+1$, where $x(P_R)$ and $x(P_L)$ are the $x$-coordinates. (see \cite[\S 2.2]{GK17})
	\item For $n=2$, we say such a polytope is \gkt{} if $\Delta$ is a triangle with vertices $(0,0)$, $P_L$ and $P_R$, with $P_L$ and $P_R$ and $(0,1)$ collinear, and $x(P_L)<0<x(P_R)<x(P_L)+1$.
	\item In both dimension $2$ and $3$, define the {\em width} of a polytope \gkt{} to be $x(P_R)-x(P_L)$.
\end{enumerate}
\end{defi}
By definition, 3-dimensional polytope $\Delta$ \gkt{} has some evident properties:
\begin{enumerate}[label={\rm(\alph*)}]
	\item The cross sections of $\Delta$ at $x=i\in \NN$ are isosceles right triangles.
	\item Projecting $\Delta\in \R^3$ \gkt{} and of width $<1$ to $xy$-plane or $xz$-plane, and then translating by the vector $(0,-1)$ will give a triangle \gkt{} with the same width.
\end{enumerate}

We first recall the following numerical criteria from \cite{GK}, \cite{GK17} for the weights for $\PP(a,b,c,d)$ or $\PP(a,b,c)$ to have a polytope \gkt{}. We rephrase the criteria as follows:
\begin{lemma} 
	\label{abcdPolytope}
	\begin{enumerate}[label={\rm(\roman*)}] 
		\item Given $w\in \Q\cap (0,1)$. Consider $\PP(a,b,c)$ with $a,b,c$ pairwise coprime. Then $\PP(a,b,c)$ has a polytope $\Delta$ \gkt{} of width $w$ if and only if there exist a relation $(e,f,-g)$ with $ae+bf=cg$ (up to a permutation of the weights $a,b,c$) and $w=cg^2/ab$. Furthermore, up to switching $a$ with $b$, and up to a shear transformation $(x,y)\mt (x,y+kx)$ for some $k\in\Z$, $\Delta$ has vertices given by (\ref{Spoly}), i.e.,
\begin{align}
		(0,0), \quad \left(-\frac{eg}{b},-\frac{er-b}{b}\right), \quad \left(\frac{fg}{a},\frac{fr+a}{a}\right),
\end{align}
where $r$ is the unique integer such that $1\leq r\leq g$, $g\mid er-b$ and $g\mid fr+a$ \cite[Prop. 5.1]{zh2017n}, and $\Delta$ is normal to the fan with ray generators given in (\ref{Sfan}).
In particular, when $w<1$, the numbers of lattice points on slices of $\Delta$ are determined by $a,b,c$.

		\item Given $W\in \Q\cap (0,1)$. Consider $\PP(a,b,c,d)$ with every $3$ weights relatively prime. Then $\PP(a,b,c,d)$ has a polytope $\Delta$ \gkt{} of width $W$ if and only if there exist positive integers $e,f,g_1,g_2$ such that up to a permutation of the weights $a,b,c$ and $d$, we have
			\[ae+bf=cg_1=dg_2, \quad W=(dg_2)^3/(abcd), \quad \gcd(e,f,g_1)=\gcd(e,f,g_2)=\gcd(g_1,g_2)=1.\]
	\end{enumerate}
\end{lemma}

The following definition is from \cite{GK17}:
\begin{defi}\label{slice}
	\cite[\S 2.2]{GK17} Suppose $\Delta$ is a $2$ or $3$-dimensional polytope \gkt{}. Suppose $m$ is a positive integer such $m\Delta$ is a lattice polytope. For any integer $i$ such that $m\cdot x(P_L)\leq i\leq m\cdot x(P_R)$, the {\em slice} at $x=i$ is the set of lattice points in $m\Delta$ with $x$-coordinates $i$. 
	When $\dim\Delta=2$, a slice of $m\Delta$ consists of consecutive lattice points on a line. When $\dim\Delta=3$, a slice of $m\Delta$ forms a right triangle with the same number $n$ of lattice points on each side. Then say the slice at $x=i$ has size $n$. 
\end{defi}

To avoid ambiguity, in the following we use $\Gamma$ to represent a $2$-dimensional polytope of \gkt{}.
We recall the following criteria in \cite{GK} and \cite{GK17} for $\bl_p X$ to be not a MDS where $X$ is a toric surface or toric 3-fold with a polytope \gkt{}.
\begin{theorem}
	{\em \cite[Thm. 1.2]{GK}} Suppose $S$ is a toric surface with fan $\Sigma$ in $\R^2$. Suppose $\Gamma\subset \R^2$ is a triangle \gkt{} with width $w$ and normal fan $\Sigma$. Let $m>0$ be a sufficiently large and divisible integer so that $m\Gamma$ is a lattice triangle. Then $\bl_p S$ is not a MDS if the following hold:
	\begin{enumerate}[label={\rm(\roman*)}]
		\item Let the slice at $m\cdot x(P_L)+1$ of $m\Gamma$ have exactly $n$ elements. Then the slice at $m\cdot x(P_R)-n+1$ of $m\Gamma$ has exactly $n$ elements.
	\item $ns_2\not\in\Z$, where $s_2:=(y(P_R)-y(P_L))/w$ is the slope of the line through $P_L$ and $P_R$.
	\end{enumerate}
	\label{gk2}
\end{theorem}
\begin{theorem}
	{\em\cite[Cor. 2.5]{GK17}} Suppose $X$ is a toric 3-fold with fan $\Sigma$ in $\R^3$. Suppose $\Delta\subset \R^3$ is a polytope \gkt{} with width $W$ and normal fan $\Sigma$. Let $m>0$ be a sufficiently large and divisible integer so that $m\Delta$ is a lattice polytope. Then $\bl_p X$ is not a MDS if the following hold:
	\begin{enumerate}[label={\rm(\roman*)}]
		\item Let the slice at $m\cdot x(P_L)+1$ of $m\Delta$ have size $n$. Then the slice at $m\cdot x(P_R)-n+1$ of $m\Delta$ has size $n$.
	\item $n(s_y,s_z)\not\in\Z^2$, where $s_y:=(y(P_R)-y(P_L))/W$ and $s_z:=(z(P_R)-z(P_L))/W$ are the $y,z$-slopes of the line through $P_L$ and $P_R$.
	\end{enumerate}
	\label{gk3}
\end{theorem}

Now a natural question is that whether there are examples of $\PP(a,b,c,d)$ meeting assumptions in Theorem \ref{main} and \cite[Cor. 2.5]{GK17}. The following proposition provides a precise answer on the overlap: 

\begin{prop}
	Suppose  $\PP(a,b,c,d)$ has a polytope $\Delta$ \gkt{} and satisfies the assumptions including (i) - (iv) of Theorem \ref{main}. Then $d=cg$, where $(e,f,-g)$ is the unique relation between $(a,b,c)$ with $w<1$.
	
	Conversely, every weighted projective 3-space $\PP(a,b,c,cg)$ such that $(a,b,c)$ has a relation $(e,f,-g)$ with $w<1$, and $\PP(a,b,c)$ has a polytope satisfying the conditions in \cite[Thm. 1.2]{GK} with width $w$, will satisfy the assumptions in both Theorem \ref{main} and \cite[Cor. 2.5]{GK17}.
	\label{comparison}
\end{prop}

\begin{remark}
	In the proof of Theorem \ref{main}, we in fact showed that weighted projective spaces $\PP(a,b,c,d)$ meeting the conditions of the theorem must contain the weighted projective plane $S=\PP(a,b,c)$ where $\bl_p S$ is not a MDS. Recall Theorem \ref{pureintersection} that $S$ is the Zariski closure of the subtorus $T_{12}=L_{12}\otimes \CC^*$, where $(L_{12})_\R$ is the $xy$-plane.

	Question: {\em Is there any $\PP(a,b,c,d)$ such that $\bl_p \PP(a,b,c,d)$ is not a MDS, but for any $2$-dimensional subtorus $T'$ of the open torus $T_N$, the blow-up $\bl_p \overline{T'}$ of the Zariski closure of $T'$ is a MDS?}

Note that the Zariski closure $\overline{T'}$ may have Picard number $1$ or $2$.
\end{remark}

\vem
We first prove Lemma \ref{abcdPolytope}. We note the following fact:
\begin{lemma}(See \cite[\S 1]{GK})
	Suppose $a,b,c$ are pairwise coprime positive integers. Then there exist at most one relation $(e,f,-g)$ of $(a,b,c)$ with $cg^2<ab$, even when permuting $a,b,c$.
	\label{uniquerelation}
\end{lemma}

\pfof{Lemma \ref{abcdPolytope}}.
First we prove (i). Suppose $\PP(a,b,c)$ has a relation of weight $w<1$, then the polytope in (\ref{Spoly}) is \gkt{} with width $w$. Conversely, suppose $S=\PP(a,b,c)$ has a polytope $\Gamma$ \gkt{} with width $w<1$. Then $S$ has a fan $\Sigma_S$ normal to $\Gamma$. Say the ray generators of $\Sigma_S$ is $r_i=(x_i,y_i)$, $i=1,2,3$. Then we can assume $y_1<0$, $y_2<0$, $y_3>0$, $ar_1+br_2+cr_3=0$, and $P_L=s(y_1,-x_1)$, $P_R=t(-y_2,x_2)$ for some $s,t\in \Q$. Since $r_i$ span the fan of $\PP(a,b,c)$, the absolute values of the $2\times 2$ minors of the following matrix should equal to $(c,b,a)$ respectively:
\[\begin{pmatrix}
		x_1& x_2 & x_3\\
		y_1& y_2 & y_3
\end{pmatrix}.\]
Now the collinearity of $P_L$, $P_R$ and $(0,1)$ gives  $w=P_R-P_L=stc$. The condition that $\overline{P_L P_R}$ being perpendicular to $r_3$ gives $bs=at=\vv{y_3}=y_3$. Therefore $w=stc=y_3^2 c/ab<1$. So $ay_1+by_2+cy_3=0$ and $y_3^2 c/ab<1$. Now $\gcd(a,b,c)=1$, so $\gcd(y_1,y_2,y_3)=1$. Write $y_1=-e, y_2=-f$ and $y_3=g$. By Lemma \ref{uniquerelation},  $(e,f,-g)$ is the unique relation. After a shear transformation of the form  $(x,y)\mt (x,y+kx)$ for some $k\in\Z$, we can assume $1\leq x_3\leq g$. Then $gx_1=ex_3\pm b$ and $gx_2=-fx_3\mp a$. So up to switching $a$ with $b$, $x_3$ is the unique integer $r$ such that  $1\leq r\leq g$, $g\mid er-b$ and $g\mid fr+a$. This shows that $\Gamma$ is of the required form, up to a reflection about the $y$-axis and a shear transformation.
The shear transformations add the same integer $k$ to the slopes of sides of $\Gamma$. Hence the numbers of lattice points on the slices are unchanged.

Next we prove (ii). Suppose $\PP(a,b,c,d)$ has a polytope $\Delta$ \gkt{}, with $P_R=(x,y,z)$, $x>0$ and $P_L=\lambda(x,y,z)$ for some $\lambda<0$. The fan $\Sigma$ is normal to $\Delta$. Therefore the four rays $R_1,\cdots, R_4$ of $\Sigma$ are the outer normal vectors of the four faces of $\Delta$. Direct calculation shows that $R_i$ is spanned by the vector $r_i$:

\begin{equation}
\begin{aligned}	
\label{rays} 
	r_1&=(1-y-z,x,x),  &\quad& r_2=(\lambda y+\lambda z-1,-\lambda x,-\lambda x),\\
	r_3&=(y,-x,0), &\quad & r_4=(-\lambda z,0,\lambda x).
\end{aligned}
\end{equation}

Now let $r'_i$ be the first lattice point in the ray $R_i$. Because $x>0$ and $\lambda<0$, there must exist positive integers $e,f,g_1,g_2$ and integers $R,S,T,U$ such that
\begin{align*}
	r'_1&=(R,e,e),\quad   r'_2=(S,f,f),\quad  r'_3 =(T,-g_1,0), \quad r'_4=(U,0,-g_2).
\end{align*}
Since $\Sigma$ is the fan of $\PP(a,b,c,d)$, up to a permutation of the weights, we have $ar'_1+br'_2+cr'_3+dr'_4=0$. Take the last two components, we have $ae+bf=cg_1=dg_2$.
Since $\Sigma$ is a fan of $\PP(a,b,c,d)$, the weights $(a,b,c,d)$ equal to the $3\times 3$ minors of the matrix with rows $r'_1,\cdots, r'_4$.
For any $3$ vectors $v_1,v_2$ and $v_3$ in $\R^3$, we denote by $\det(v_1,v_2,v_3)$ the determinant of the square matrix with row vectors $v_1,v_2$ and $v_3$.
Then we have
\begin{align*}
	a&=\vv{\det(r'_2,r'_3,r'_4)}=\frac{g_1 g_2}{x}\vv{Sx+fy+fz}=\frac{g_1 g_2}{x}\left|\frac{(\lambda y+\lambda z-1)f}{-\lambda}+fy+fz\right|=-\frac{fg_1 g_2}{\lambda x},\\
	b&=\vv{\det(r'_1,r'_3,r'_4)}=\frac{g_1 g_2}{x}\vv{Rx+ey+ez}=\frac{g_1 g_2}{x}\vv{(1-y-z)e+ey+ez}=\frac{eg_1 g_2}{x},
\end{align*}
where we used that each $r'_i$ is a scalar multiple of $r_i$. Note that the other two equations of $c$ and $d$ do not give new algebraic relations. As a result,
\begin{gather}
	x=\frac{eg_1 g_2}{b}, \quad \lambda=-\frac{bf}{ae},\label{alambda}\\
	W=x(P_R)-x(P_L)=x-\lambda x=\frac{eg_1 g_2}{b}\left(1+\frac{bf}{ae}\right)=\frac{eg_1 g_2}{b}\cdot \frac{dg_2}{ae}=\frac{cg_1 \cdot dg_2\cdot dg_2}{abcd}=\frac{(dg_2)^3}{abcd}.\label{Wform}
\end{gather}

At last, the coprime conditions follow from the assumption that every $3$ of $a,b,c,d$ are relatively prime, and the expression of $a,b,c,d$ as the determinants of $r'_i$ with $R,S,T$ and $U$ are integers.
This proves the `only if' direction. Conversely, suppose $ae+bf=cg_1=dg_2$ and $W=(dg_2)^3/(abcd)$. We can always choose integers $T$ and $U$ such that $\gcd(T,g_1)=\gcd(U,g_2)=1$. Let
$y=Tx/g_1$ and $z=Ux/g_2$, with $x$ and $\lambda$ given above in (\ref{alambda}).  The parameters $x,y,z,\lambda$ determine a fan $\Sigma'$ with rays $r_i$ from (\ref{rays}), and a polytope $\Delta'$ with  $P_R=(x,y,z)$, $x>0$ and $P_L=\lambda(x,y,z)$. Then it is straightforward that $\Sigma'$ is a fan of $\PP(a,b,c,d)$, and $\Delta'$ is \gkt{} with width $W$, whose normal fan is $\Sigma'$. This proves the `if' direction. 
\qed

Finally we prove Proposition \ref{comparison}.

\pfof{Proposition \ref{comparison}.}  Suppose $\PP(a,b,c,d)$ has a polytope $\Delta$ \gkt{} and meets the assumptions of Theorem \ref{main}. Then by Lemma \ref{abcdPolytope}, there exist $e,f,g_1,g_2\in \Z_{>0}$ such that $ae+bf=cg_1=dg_2$ (up to a permutation of the weights $a,b,c$ and $d$), and the width $W$ of $\Delta$ equals $(dg_2)^3/(abcd)\leq 1$.
In this equation, $a$ and $b$ are symmetric. The weights $c$ and $d$ are also symmetric. Hence up to symmetry either $\bl_p \PP(a,b,c)$ is not a MDS or $\bl_p \PP(b,c,d)$ is not a MDS.

Case I. $\bl_p \PP(a,b,c)$ is not a MDS, with relations $(E,F,-G)$ such that the width $w<1$. By the argument above,
\[1\geq W=\frac{(dg_2)^3}{abcd}=\frac{cg_1^2 g_2}{ab}.\]
We claim $W<1$. Otherwise $W=1$. Then $cg_1^2 g_2=ab$, so $c\mid ab$, which contradicts the assumption of Theorem \ref{general} that $a,b,c$ are pairwise coprime. 

Hence $cg_1^2/ab<1/g_2\leq 1$. By Lemma \ref{abcdPolytope}, $\gcd(e,f,g_1)=1$. Now $(e,f,-g_1)$ is a relation between $(a,b,c)$ with $\gcd(e,f,-g_1)=1$ and width $c(g_1)^2/(ab)=cg_1^2/(ab)<1$.  By Lemma \ref{uniquerelation}, we must have $e=E,f=F$ and $g_1=G$, $ae+bf=cg_1$, and the width of $(e,f,-g_1)$ is
\[w=\frac{cG^2}{ab}=\frac{cg_1^2}{ab}<\frac{1}{g_2}\leq 1.\]
Suppose $g_2\geq 2$. Then $w\leq 1/2$. By Theorem 2.5 and 2.6 of \cite{zh2017n}, if $w\leq 1/2$, then $\bl_p \PP(a,b,c)$ is a MDS, which contradicts the assumption. Therefore $g_2=1$, and $d=cg_1$. 

Case II. $\bl_p \PP(b,c,d)$ is not a MDS, and $\gcd(b,c,d)=1$. This together with $cg_1=dg_2$ implies that $g_1=kd$ and $g_2=kc$ for some $k\in \Z_{>0}$. Now
\[1\geq W=\frac{cg_1^2 g_2}{ab}=\frac{k^3 c^2 d^2}{ab}.\]
Hence $k^3 c^2 d^2\leq ab$. On the other hand, $kcd=cg_1=ae+bf\geq a+b\geq 2\sqrt{ab}$. Hence $k^3 c^2 d^2\geq k\cdot (4ab)>ab$, so we reached a contradiction. This shows Case II does not happen and proves the first half of Proposition \ref{comparison}.

Next we prove the second half of Proposition \ref{comparison}.
Consider any $S=\PP(a,b,c)$ such that $a,b,c$ are pairwise coprime, $(e,f,-g)$ is a relation between $(a,b,c)$ of width $w<1$ and $S$ satisfies the assumptions in \cite[Thm. 1.2]{GK}. Then $\bl_p \PP(a,b,c)$ is not a MDS.

Now $X:=\PP(a,b,c,cg)$ satisfies conditions (i), (ii) and (iv) of Theorem \ref{main}. Since $d=cg$, we have $d^2 w/(abc)=cg^2 w/(ab)=w^2<1$. This verifies condition (iii). Hence $X=\PP(a,b,c,cg)$ is an example of Theorem \ref{main}.

It remains to show that $X=\PP(a,b,c,cg)$ satisfies the two assumptions in \cite[Cor. 2.5]{GK17}. Indeed, here $ae+bf=cg=d\cdot 1$ with $cg^2/ab<1$. By Lemma \ref{abcdPolytope}, $X$ and $S=\PP(a,b,c)$ have polytopes $\Delta$ and $\Gamma$ \gkt{}. Let $r$ be the unique integer such that $1\leq r\leq g$, $g\mid er-b$ and $g\mid fr+a$.  Recall the proof of Lemma \ref{abcdPolytope}. By setting $T=-r$ and $U=0$, we can determine the parameters $x,y,z$ and $\lambda$ to give
\[P_L=\left(-\frac{fg}{a},\frac{fr}{a},0\right), \quad P_R=\left(\frac{eg}{b},-\frac{er}{b},0\right).\]
This gives a polytope $\Delta$ \gkt{}.
The fan $\Sigma$ of $X$ can be chosen as the fan with ray generators
\[	r'_1=\left(\frac{er-b}{g},e,e\right),\quad r'_2=\left(\frac{fr+a}{g},f,f\right),\quad r'_3=(-r,-g,0),\quad r'_4=(0,0,-1).\]
Define $\Gamma$ to be the projection of $\Delta$ to the $xy$-plane, after translating $(0,1)$ to $(0,0)$ and a reflection about $y$-axis. Then $\Gamma$ is the triangle given by (\ref{Spoly}), which is a polytope of $S=\PP(a,b,c)$. 

Now let $\Gamma'$ be the reflection of $\Gamma$ about the $y$-axis. By the hypothesis and Lemma \ref{abcdPolytope} (i), either $(S,\Gamma)$ or $(S, \Gamma')$ meets the assumptions of \cite[Thm. 1.2]{GK}. By symmetry we can assume the case $(S, \Gamma)$. Then \cite[Thm. 1.2]{GK} (i) says that for some $m>0$, the slice at $m\cdot x(P_L)+1$ of $m\Gamma$ has exactly $n$ elements, and the slice at $m\cdot x(P_R)-n+1$ of $m\Gamma$ has exactly $n$ elements too. By Definition \ref{slice}, every slice of $\Delta$ forms a right triangle with the same number of lattice points on each right side. Hence, both slices of $m\Delta$ at $m\cdot x(P_L)+1$ and $m\cdot x(P_R)-n+1$ of $m\Delta$ have size $n$. This shows that (i) of  \cite[Cor. 2.5]{GK17} holds. For (ii) of \cite[Cor. 2.5]{GK17}, we have $s_y$ equals $s_2$ of the triangle $\Gamma$ in $xy$-plane. 
If $\Gamma$ meets the assumption (ii) of \cite[Thm. 1.2]{GK}, then $ns_y=ns_2\not\in\Z$, so $\Delta$ meets the assumption (ii) of \cite[Cor. 2.5]{GK17}. Therefore, $X$ satisfies the two assumptions in \cite[Cor. 2.5]{GK17}.\qed

\begin{remark}
	Consider $X=\PP(a,b,c,cg)$ in the overlap described in Proposition \ref{comparison}. A comparison with \cite[Lem. 5.1, 5.2]{GK17} shows that the curve $C\subset \bl_p X$ we constructed in Definition \ref{CC0}, whose class is extremal in the Mori cone $\NEb(\bl_p X)$ (by Theorem \ref{general}), is the same curve $C$ constructed in \cite[Lem. 5.1, 5.2]{GK17}.  
\end{remark}

\begin{example}
	An example in such family of $\PP(a,b,c,cg)$ is $\PP(7,15,26,52)$. By \cite{GK}, $\bl_p \PP(7,15,26)$ is not a MDS. The relation is $(e,f,-g)=(1,3,-2)$. Both Theorem \ref{main} and \cite[Cor. 2.5]{GK17} apply to $\PP(7,15,26,52)$, so $\bl_p \PP(7,15,26,52)$ is not a MDS. 
\end{example}

\vem

%% file: application.tex
\section{Application}\label{extra}

We apply Proposition \ref{g} to the following examples in \cite{AGK17}.
By \cite[Ex. 1.4]{AGK17}, the blow-up $\bl_p S$ of the following $S=\PP(a,b,c)$ at the identity point $p$ is not a MDS:
\begin{align}
	(a,b,c)=((m+2)^2,(m+2)^3+1,(m+2)^3 (m^2+2m-1)+m^2+3m+1),
	\label{higherM}
\end{align}
where $m$ is a positive integer.

We briefly review the geometry on those $\bl_p S$.
By \cite[Thm. 1.1]{AGK17}, for every positive integer $m\geq 1$, there exists an irreducible polynomial $\xi_m\in \CC[x,y]$ such that $\xi_m$ has vanishing order $m$ at $(1,1)$ and the Newton polygon of $\xi_m$ is a triangle with vertices $(0,0), (m-1,0)$ and $(m,m+1)$.
Now the weighted projective plane $S$ above satisfies the conditions of \cite[Thm. 1.3]{AGK17}. Then by \cite[Thm. 1.3]{AGK17} and its proof, the polynomial $\xi_m$ above defines a curve $H$ in $S$, passing through $p$ with multiplicity $m$, such that the proper transform $C$ of $H$ in $\bl_p S$ is a negative curve. Then $C\neq e$.
The proof of \cite[Thm. 1.3]{AGK17} in fact shows that $H$ is the polarization given by the triangle $\Delta$ with vertices $(-\alpha,0), (m-1+\beta,0),(m,m+1)$, with
\[\alpha=\frac{1}{(m+2)^2}, \quad\beta=\frac{(m+2)^2+1}{(m+2)^3+1}.\] 
Therefore on $S$ we have 
\[H^2=2\mathrm{Area}(\Delta)=\frac{(m+1)^2 c}{ab}.\]
Let $B$ be the pseudo-effective divisor on $S$ generating $\Cl(S)\cong\Z$. Then $H\sim rB$ for some $r\in \Q_{>0}$. Since $B^2=1/abc$ and $H^2=r^2 B^2$, we have $r=c(m+1)$, so $[H]= c(m+1)[B]\in \Cl(S)$. Therefore $C\sim c(m+1)\pi^*B-me$.

When $m\geq 2$, those $S$ above have width $w\geq 1$, so Theorem \ref{main} does not apply to $S$. Nevertheless, by Proposition \ref{g}, we have the following examples:
\begin{corollary}\label{mgeq1}
	Let $X=\PP(a,b,c,d_1,d_2,\cdots,d_{n-2})$ where
\[(a,b,c)=((m+2)^2,(m+2)^3+1,(m+2)^3 (m^2+2m-1)+m^2+3m+1),\]
such that $m\in \Z_{>0}$, every $d_i$ lies in the semigroup generated by $a,b$ and $c$, and that every \linebreak $d_i< abm/(m+1)$. Let $p$ be the identity point of the open torus in $X$. Then $\bl_p X$ is not a MDS.
\end{corollary}

%% file: MDSBlowupWPS.bbl
\begin{thebibliography}{GAGK17}

\bibitem[BCHM10]{BCHM}
Caucher Birkar, Paolo Cascini, Christopher~D. Hacon, and James McKernan.
\newblock Existence of minimal models for varieties of log general type.
\newblock {\em Journal of the American Mathematical Society}, 23(2):405--468,
  2010.

\bibitem[Cas15]{castravet16}
Ana-Maria Castravet.
\newblock Mori {D}ream {S}paces and blow-ups.
\newblock {\em Proceedings of the AMS Summer Institute in Algebraic Geometry
  2015}, pages 143--168, 2015.

\bibitem[CLS11]{coxtoric}
David~A. Cox, John~B. Little, and Henry~K. Schenck.
\newblock {\em Toric Varieties}.
\newblock Graduate studies in mathematics. American Mathematical Society, 2011.

\bibitem[CT15]{castravettevelev2015}
Ana-Maria Castravet and Jenia Tevelev.
\newblock $\overline {M}_{0,n}$ is not a {Mori} dream space.
\newblock {\em Duke Mathematical Journal}, 164(8):1641--1667, 06 2015.

\bibitem[Cut91]{cutkosky1991}
Steven~Dale Cutkosky.
\newblock Symbolic algebras of monomial primes.
\newblock {\em J. Reine Angew. Math}, 416:71--89, 1991.

\bibitem[FS97]{FS1997}
William Fulton and Bernd Sturmfels.
\newblock Intersection theory on toric varieties.
\newblock {\em Topology}, 36(2):335--353, 1997.

\bibitem[Ful93]{introTV}
William Fulton.
\newblock {\em Introduction to Toric Varieties}.
\newblock Princeton University Press, 1993.

\bibitem[GAGK17]{AGK17}
Javier Gonz\'{a}lez-Anaya, Jos\'{e}~Luis Gonz\'{a}lez, and Kalle Karu.
\newblock On a family of negative curves, 2017,
  \href{https://arxiv.org/abs/1712.04635}{arxiv:1712.04635}.

\bibitem[GK16]{GK}
Jos\'{e}~Luis Gonz\'{a}lez and Kalle Karu.
\newblock Some non-finitely generated {Cox} rings.
\newblock {\em Compositio Mathematica}, 152(5):984--996, 2016.

\bibitem[GK17]{GK17}
Jos\'{e}~Luis Gonz\'{a}lez and Kalle Karu.
\newblock Examples of non-finitely generated {Cox} rings, 2017,
  \href{https://arxiv.org/abs/1708.09064}{arxiv:1708.09064}.

\bibitem[GNW94]{gnw1994}
Shiro Goto, Koji Nishida, and Kei-ichi Watanabe.
\newblock Non-{C}ohen-{M}acaulay symbolic blow-ups for space monomial curves
  and counterexamples to {C}owsik's question.
\newblock {\em Proceedings of the American Mathematical Society},
  120(2):383--392, 1994.

\bibitem[He17]{zh2017n}
Zhuang He.
\newblock {New examples and non-examples of MDS when blowing up toric
  surfaces}, 2017, \href{https://arxiv.org/abs/1703.00819}{arxiv:1703.00819}.

\bibitem[HK00]{hu2000}
Yi~Hu and Sean Keel.
\newblock Mori dream spaces and {GIT}.
\newblock {\em The Michigan Mathematical Journal}, 48(1):331--348, 2000.

\bibitem[HKL16]{hkl}
J\"urgen Hausen, Simon Keicher, and Antonio Laface.
\newblock On blowing up the weighted projective plane, 2016,
  \href{https://arxiv.org/abs/1608.04542}{arXiv:1608.04542 [math.AG]}.

\bibitem[KM08]{kollar2008}
Janos Koll{\'a}r and Shigefumi Mori.
\newblock {\em Birational Geometry of Algebraic Varieties}.
\newblock Cambridge University Press, 2008.

\bibitem[Muk05]{Mukai05}
Shigeru Mukai.
\newblock Finite generation of the {N}agata invariant rings in {A-D-E} cases.
\newblock In {\em RIMS Preprint 1502}, 2005.

\bibitem[Sri91]{HS}
Hema Srinivasan.
\newblock On finite generation of symbolic algebras of monomial primes.
\newblock {\em Communications in Algebra}, 19(9):2557--2564, 1991.

\end{thebibliography}
